# Paddy: Evolutionary Optimization Algorithm for Chemical Systems and Spaces

Armen Beck[1], Jonathan Fine[1], Gaurav Chopra[1,2,3,4,5,6,]*

[1]Department of Chemistry and Computer Science (*by courtesy*), Purdue University, 720 Clinic Drive, West Lafayette, IN 47907

[2]Purdue Institute for Drug Discovery, West Lafayette, IN 47907

[3]Purdue Center for Cancer Research, West Lafayette, IN 47907

[4]Purdue Institute for Inflammation, Immunology and Infectious Disease, West Lafayette, IN 47907

[5]Purdue Institute for Integrative Neuroscience, West Lafayette, IN 47907

[6]Regenstrief Center for Healthcare Engineering, West Lafayette, IN 47907

*Corresponding author email – gchopra@purdue.edu

## Abstract

Optimization of chemical systems and processes have been enhanced and enabled by the guidance of algorithms and analytical approaches.  While many methods will systematically investigate how underlying variables govern a given outcome, there is often a substantial number of experiments needed to accurately model these relations.  As chemical systems increase in complexity, inexhaustive processes must propose experiments that efficiently optimize the underlying objective, while ideally avoiding convergence on unsatisfactory local minima.  We have developed the Paddy software package around the Paddy Field Algorithm, a biologically inspired evolutionary optimization algorithm that propagates parameters without direct inference of the underlying objective function.  Benchmarked against the Tree of Parzen Estimator, a Bayesian algorithm implemented in the Hyperopt software Library, Paddy displays efficient optimization with lower runtime, and avoidance of early convergence.  Herein we report these findings for the cases of: global optimization of a two-dimensional bimodal distribution, interpolation of an irregular sinusoidal function, hyperparameter optimization of an artificial neural network tasked with classification of solvent for reaction components, and targeted molecule generation via optimization of input vectors for a decoder network.  We anticipate that the facile nature of Paddy will serve to aid in automated experimentation, where minimization of investigative trials and or diversity of suitable solutions is of high priority.

## Introduction

Optimization is used in all of chemical sciences, including identifying synthetic methodology[1–3], chromatography[4–6] conditions, calculating transition state geometry[7], or selecting materials and drug formulations[8–11]. Typically, several parameters or variables need to be optimized that are done either by human chemists using chemical intuition or computational methods to identify suitable conditions[12–14]. The development of automated optimization procedures for repetitive human tasks in chemical sciences, such as shimming,[15] chromatograph peak assignment[16] and developing bioanalytical workflows,[17] have saved time and resources. Several chemical optimization methods have been used iteratively in a task-specific manner to optimize an objective to model or select experimental conditions for chemical and biological processes[18–28]. However, iterations with stochastic optimization algorithms have been shown to provide a better alternative to deterministic algorithms to find optimal solutions[29–32]. Examples includes the use of stochastic gradient decent algorithm that outperforms gradient decent to find optimal solutions[33-35].

Several artificial intelligence and machine learning (AI/ML) architectures have been used in the chemical sciences where stochastic optimization algorithms are needed for train-validate-test cycles[36,37]. These AI/ML algorithms are used in several areas of the chemical sciences, such as retrosynthesis[38], reaction condition prediction[39–42], catalyst design[43,44], drug design[45–48], spectral interpretation[49–51], retention time prediction[52], and for molecular simulations[53–55]. Specific generative AI neural network architectures have also been used for inverse design[56,57] and property-specific generation of molecules[45,58–60]. For optimization tasks related to laboratory automation, that use closed-loop procedures, several methods have been developed[61] including active learning using neural networks[62,63]. In addition, the use of Bayesian methods[64,65], genetic algorithms[66] and other iterative optimization methods[67] have resulted in useful chemical outcomes without the use of prior learning.

Evolutionary algorithms are a class of optimization methods inspired by biological evolution using a starting set of possible solutions (seeds) to the problem that are then evaluated using a 'fitness (objective) function' to 'evolve' the next set of solutions iteratively towards identifying optimal solutions. Using directed sampling to maximize a fitness function, the evolutionary optimization algorithms propagate parameters to find the set of optimal solutions for a given problem. Several types of evolutionary algorithms include, the most popular genetic algorithms, evolution strategies, differential evolution, and estimation of distribution algorithms[34]. The propagation between iterations use (meta)heuristic approaches with a set of rules that include simulated annealing[68], genetic algorithms[69], Tabu search[70], hill climbing methods[71], and particle swarm[72] to name a few. In contrast, Bayesian methods lend to directed optimization, guided by sequential updates of a probabilistic model and inferring the return on sampling, often via an acquisition function[73]. Furthermore, Bayesian optimization methods have also been reported in the chemical literature for the optimization of neural networks[41], generative sampling[74,75], and as a general-purpose optimizer for chemistry[64,76,77].

Herein we have implemented a new class of evolutionary algorithm, the Paddy field algorithm (PFA)[78] as a Python library, named Paddy, which includes heuristic methods that operate on a reproductive principle dependent on solution fitness and the distribution of population density among a set of selected solutions. We also show the advantages of using Paddy, when compared to Bayesian optimization[79] implemented in the Hyperopt library[80] with random solutions as controls. We compared test cases for accuracy, speed, sampling parameters and sampling performance across various optimization problems. Specifically, the problems include identification of the global maxima of a two-dimensional bimodal distribution, interpolation of an irregular sinusoidal function, hyperparameters optimization of a neural network trained on chemical reaction data, and comparison of performance for targeted molecule generation using junction-tree variational autoencoder. Paddy outperforms the Bayesian implementation and resulted in effective sampling of conditions to identify optimal solutions. Furthermore, Paddy was designed with user experience in mind with the features to save and recover Paddy trials. We include a complete documentation and code via GitHub (https://github.com/chopralab/paddy) to encourage others to use and extend Paddy for their chemical optimization tasks.

## Methods

### Formulation of the Paddy Field Algorithm (PFA)

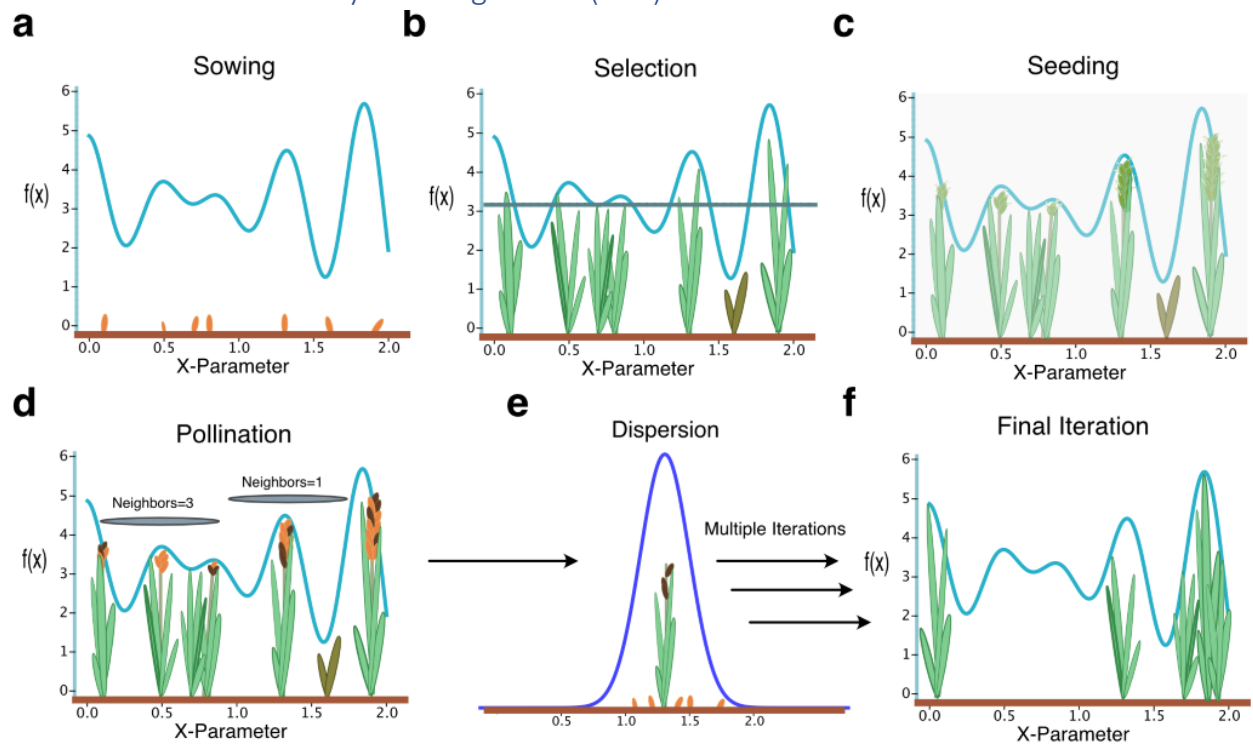

**Figure 1**. **Overview of the Paddy Field Algorithm.** Paddy is initiated by the sowing step (a) where objective function parameters, arbitrary in dimensionality, are randomly sown as the initial population of seeds. After evaluation of the seeds, the selection step (b) applies a selection operator to select a, user defined,

number of top preforming plants to further propagate. The seeding step (c) then calculates how many seeds a selected plant should respectively generate as to account for fitness across parameter space, such as fertility of soil determines the number of flowers a plant can grow. The pollination step (d) then reinforces the density of selected plants by eliminating seeds proportionally for those with fewer than the maximum number of neighboring plants within Euclidian space of the objective function variables. The sowing step (e) then assigns new parameter values to pollinated seeds by randomly dispersing across a Gaussian distribution, with the mean being the parameter values of the parent plant. The algorithm terminates (f) after converging or running for the number of iterations set.

The PFA was inspired by the reproductive behavior of plants that is based on the relationship of soil quality, pollination and plant propagation to maximize plant fitness. The PFA is developed without knowing this underlying relationship to iteratively optimize a fitness (objective) function using a five-phase process (a-e, **Figure 1**). First, for any objective (fitness) function, $y = f(x)$, with dependent parameters $(x)$ of n-dimensions, PFA treats individual parameters $x = \{x_1, x_2, \ldots, x_n\}$ as seeds to define a numerical propagation space. Next, these seeds are converted to plants by evaluating the objective (fitness) function, $y = f(x)$, at the respective seed values. The resulting evaluation provides plant fitness score values thereby assessing soil quality. Parameters $(x_H \in x)$ that result in plants of high fitness $(y_H \in y)$ are further evaluated and selected for seeding and propagation $(y^* \in y_H)$. The number of neighboring plants and their fitness scores determine the number of seeds in each round $(s)$ produced by a plant selected for propagation $(y^* \in y_H)$ thereby directing plant density mediated pollination. The parameter values $(x^* \in x)$ for selected plants are then modified by sampling from a Gaussian distribution. We provide the details of the five-phase process (a-e) as follows:

a. **Sowing**: The Paddy algorithm is initiated with a random set of user defined parameters $(x)$, as starting seeds for evaluation.

b. **Selection**: The fitness function, $y = f(x)$, is evaluated for the selected set of seed parameters $(x)$ to convert seeds to plants. A user-defined threshold parameter $(H)$ that defines the selection operator which selects the number of plants based on the sorted list of evaluations $(y_H)$ for respective seeds $(x_H)$. These function evaluations can also be taken from previous iterations, for further propagation (eq. 1).

$$f(x) = y = \{y_{min}, \ldots, y_{max}\},$$
$$H[y] = H[f(x)] = f(x_H) = y_H = \{y_t, \ldots, y_{max}\} \: \forall \: x_H \in x, \: y_H \in y \ldots \ldots \ldots \ldots (1)$$

where $y_H$ is the sorted list of function evaluations (selected plants) from all current and previous evaluations $y$ satisfying the threshold $H$ for the set of seeds or parameters $x_H$ that belong to all parameters $x$.

c. **Seeding**: The plants are further selected $y^* \in y_H$ to calculate the number of potential seeds ($s$) for propagation as the fraction of user-defined maximum number of seeds ($s_{max}$) based on the min-max normalized fitness value (eq. 2).

$$s = s_{max} \left([y^* - y_t] / [y_{max} - y_t]\right) \; \forall \; y^* \in y_H \quad \ldots \ldots \ldots \ldots \ldots \ldots \ldots \ldots \ldots (2)$$

where $s$ is the number of seeds for selected plants (function evaluation) $y^*$ that belongs to the sorted ($y_t$ minimum to $y_{max}$ maximum) list of plants satisfying the threshold $y_H$.

d. **Pollination**: This step is related to clustering based on density of all selected plants $y^* \in y_H$ (function evaluation) such that the number of seeds to be dispersed by plants (new parameters $x$ to be evaluated) is based on the number of neighbors to $y^*$. The number of neighbors, $v$, is used to calculate the pollination term $U$ (eq. 3) that ranges from 0.368 (inverse of Euler's number, $e^{-1}$) to 1 ($e^0$). The total number of pollinated seeds, $S$, to be subsequently propagated is the product of pollination term $U$ and the number of seeds for selected plants $s$ (eq. 4).

$$U = e^{\left[\frac{v}{v_{max}} - 1\right]} \quad \ldots \ldots \ldots \ldots \ldots \ldots \ldots \ldots \ldots (3)$$

$$S = U \times s \quad \ldots \ldots \ldots \ldots \ldots \ldots \ldots \ldots \ldots (4)$$

where the number of neighbors $v$ (eq. 5) is defined as the number of selected plants or function evaluations $y^* = f(x_k) \in y_H$ at $x_k \in x_H$ within the radius ($r$) of the plant or function evaluation $f(x_j)$ being considered at $x_j \ni x_k \neq x_j \; \forall \; x_i, x_k \in x_H$. When the absolute distance between plants (function evaluations) is less than the user defined hyperparameter $r$, they are all considered as neighbors (eq. 5). To this end, the number of neighbors affects pollination term $U$ in eq. 4 where term $U = 1$ for maximum number of neighbors $v_{max}$ and reduced to 0.368 (inverse of Euler's number, $e^{-1}$) for no neighbors, $v = 0$.

$$v = |n|, \; n = \{x_k \in x_H \; | \; \|x_j - x_k\| - r < 0, \; y^* = f(x_k) \in y_H, x_k \neq x_j\} \ldots \ldots \ldots \ldots \ldots \ldots (5)$$

e. **Dispersion**: For each plant (function evaluation) with pollinated seeds, the parameter values for new seeds are initialized by sampling a Gaussian distribution where the parameter values of the parent plant define the mean of the Gaussian distribution for each parameter. The standard deviation ($\sigma$) is a hyperparameter that affects the dispersion of seeds (conditions) around each selected plant.

The steps **a-e** are then repeated until the desired number of iterations or specific termination conditions are met.

## Implementation and Extension of the Paddy Field Algorithm

We have implemented the PFA algorithm and extended it with new features for chemical optimization problems. We have modified PFA[78] where the threshold parameter ($H$) is adjusted based on the user defined random seeds during initiation. This allows for maximum flexibility in selecting seeds and threshold to allow for cases when the number of random seeds are lower than the threshold during initiation of PFA. Specifically, during initiation if the number of seeds is lower than threshold number to select the seeds for the next round, the value of $H$ is equal to the rounded whole number of 75% of the number of random seeds defined by the user. In addition, the neighborhood function is modified in the pollination phase to mitigate early termination of the algorithm. For Paddy, we use Euclidean distance to determine the spatial distance between plants. The neighborhood function is dependent on the radius parameter that can result in early termination of the algorithm, in that the plants produce zero new seeds due to the radius resulting in zero neighbors. To this end, we have formulated the neighborhood function with an adaptive radius to mitigate early termination. If the initial evaluation calculates zero neighbors, the 0.75th quantile for the distance between plants is used as the radius parameter. If the 0.75th quantile radius results in zero neighbors being assigned, the quantile value is iteratively decreased by 0.05 until a nonzero number of neighbors is assigned to a plant. If the 0.05th quantile fails to generate neighbors, each plant is evaluated as having one neighbor, effectively dropping the pollination term for the given iteration.

The termination condition is defined for equal values of $y_t$ and $y_{max}$. Additionally, in Paddy, the standard deviation parameter used for the dispersion phase is defined as 0.2. To provide flexibility to the user, modifications to the algorithm have been introduced in Paddy that facilitate alternative dispersion behavior, in addition to an alternative formulation of the selection phase which are described in the subsequent section.

We have introduced several alternative methodologies to provide users greater flexibility to control different features of the algorithm that include:

- *Population Mode*: The selection phase is as described for the native PFA, where plants generated during any previous iteration are considered. As mentioned previously, population mode differs from the native PFA by having a flexible threshold parameter during, and only during, random initiation. The originally defined threshold parameter is recovered after the first iteration and remains static sense the full population of plants remain available for propagation. If the selected threshold parameter is too large compared to the number of random seeds defined during initiation, population mode may not complete. This can result as the threshold parameter will not auto-scale for low numbers of plants post random initiation.
- *Generational Mode*: The selection phase is modified such that only plants generated by the previous iteration are considered, rather than applying the threshold operator across all plants evaluated. A flexible threshold parameter is implemented as previously described, as some iterations may yield a number of seeds lower than the operator. The originally defined threshold operator is recovered and otherwise used each iteration.
- *Scaled Gaussian*: The standard deviation for the Gaussian applied during dispersion is calculated with an inherited scaling term ($\delta$) (eq. 6). The scaling term is initiated as zero and inherited in a variative manner where new values are generated by selecting from a Gaussian distribution, where the mean is the current scaling term and the standard deviation being 0.2.

$$\sigma = (0.2^{10})^{\delta} \dots\dots\dots\dots\dots\dots\dots\dots\dots\dots (6)$$

- *Parameter Type*: The parameter type determines the handling of values generated by Paddy where parameter types are either a continuous value or an integer value that is rounded after being generated.
- *Parameter Limits*: The explicit bounding of a parameter value is supported by Paddy. Limits can be either one-sided or two-sided. If parameter values are generated outside set limits, they are clamped to the limit value.
- *Parameter Normalization*: Parameters with two-sided limits can be normalized during the dispersion phase via min-max normalization with limit values.

## Min/Max Optimization of a Two-Dimensional Bimodal Distribution

The Paddy evolutionary algorithm and Hyperopt Bayesian algorithm were used to find the maxima for bimodal function with two parameters (x, y). Each algorithm was run 100 times with random initial seeds to test for robustness of the results. Hyperopt was run using the Tree-structured Parzen Estimator for 500 evaluations, and changed (x, y) parameters using 'hp.uniform' to propagate values between 0 and 1. Paddy was run in Generational mode with scaled Gaussian type setting and each (x, y) parameter limits of 0 and 1 that was randomly propagated with 0.01 resolution within the limits. PFARunner parameters were set where: the number of random seeds as 50, yt as 50, Qmax as 100, r as 0.02, and iterations being 5. Solutions for locating the global maxima were defined as values greater than 0.81 when evaluating parameters with the analytic.

## Gramacy & Lee Interpolation

Paddy, Hyperopt, and the random search algorithm were run in the same manner as for min/max optimization regarding environment, the number of executions, and random seeds. Interpolation of the Gramacy & Lee function was done using a 32$^{nd}$ degree trigonometric polynomial with 65 coefficients values ranging between –1 and 1. Interpolative performance was evaluated by calculating the mean squared error between the Gramacy & Lee function and generated trigonometric polynomial between –0.5 and 2.5 with a resolution of 0.001. The random sampling algorithm was used to generate the 65 coefficients using the Numpy

'random.uniform' function, with 5000 evaluations per execution. Hyperopt was run using the Tree-structured Parzen Estimator for 1500 evaluations, and optimized the 65 coefficients using the 'hp.uniform' to propagate values between –1 and 1. Paddy was run in Generational mode with the Gaussian type set to default and with limits of –1 and 1, and randomly propagated in range of the limits with a resolution of 0.05. PFARunner parameters were set where: the number of random seeds as 25, yt as 25, Qmax as 25, r as 0.02, and iterations being 10.

## Multilayer perceptron (MLP) Hyperparameter Optimization

Training data was obtained from the Daniel Lowe[81] repository (https://figshare.com/articles/dataset/Chemical_reactions_from_US_patents_1976-Sep2016_/5104873) and preprocessed. Briefly, the initial set of reaction SMILES, from which the subset used in this work is from, was generated by initially removing atom mapping, selecting reaction strings containing solely solvents as agents, and by associating ionic compounds with pseudo covalent bonds. Additionally, reactions where more than four reagents, post condensing of ionic pairs, or more than one product were removed. The subset used contains 4994 reactions with 30 types of solvents, and were converted into bitvectors after separating the reaction components from their respective solvent. The conversion to bitvectors was done using RDKit's 'GetMorganFingerprintAsBitVec' method[82] to produce 2048 length Morgan Fingerprints[83] using an atom radius of 2. The bitvectors and solvent labels were then converted into arrays, using onehot encoding for solvent.

The machine learning was done using the Keras package for generating and training the neural networks, while the Scikit Learn library was used for data splitting and performance assessment[84,85]. The multilayer perceptrons generated were comprised of two dense hidden layers with dropout terms and using rectified linear units as the activation function. Input and output layers were 2048 and 30 neurons in length respectively using softmax activation, and the model was compiled to use categorical crossentropy as the loss function while using the Adam optimizer. Stratified K-fold validation was used of three-fold splitting of training and validation data. Models were trained for five epochs, using batch sizes of 1000, with validation scores being calculated as micro F1 scores. The F1 scores of the three resulting models post three-fold cross validation were averaged as to provide a single value for the algorithms to optimize.

The hyperparameters of the two dropout terms, ranging from 0 – 1, and the lengths of the hidden layers, 300 – 3000 and 32 – 2000 neurons for the first and second layers respectively, were each optimized over 100 trials using Paddy, Hyperopt, and a random search algorithm. The random search algorithm generated random dropout terms and layer lengths using the Numpy 'random.uniform' and 'random.randint' functions, to propagate values for the dropout and layer length terms within their appropriate ranges, for 200 evaluations. Hypeopt was run with the Tree-structured Parzen Estimator for 200 evaluations, and used the 'hp.uniform' and 'hp.quniform' functions for dropout and layer length value generation. Paddy was run in

Generational mode with the Gaussian type set to default and using normalization when generating parameter values. Random propagation was done with dropout values between 0 – 0.5 with 0.05 in resolution and layer lengths of 300 – 3000 and 32-500 with resolutions of 0.05 for the first and second hidden layers, and propagated within the parameter limits for subsequent iterations. PFARunner parameters were set where: the random seed number as 25, yt as 5, Qmax as 10, r as 0.2, and iterations being 7.

### Junction Tree Variational Autoencoder (JTVAE) Latent Space Sampling

The JTVAE pretrained models are available on GitHub (https://github.com/wengong-jin/icml18-jtnn) that includes the conda environment and all Python v2 dependencies. With set random seeds to ensure reproducibility, latent vectors were decoded to generate SMILES strings prior to evaluation. Fitness calculated using solely Tversky Similarity, was done using the RDKit library, by converting SMILES to RDKit mol structures and subsequently to Morgan Fingerprints. The Morgan Fingerprints were generated with a bit radius of two, and length of $2^{23}$ as to minimize bit collision incidents, via the 'GetMorganFingerprintAsBitVect' method. The Morgan Fingerprints of the generated SMILES and Pazopanib were then compared via Tversky Similarity with coefficients α = 0.5 and β = 0.01. The α and β coefficients were used to weight the relative complements of the Pazopanib fingerprint in the generated fingerprint and vise versa respectively. For the random sampling algorithm, it generated the tree and graph latent vectors as two arrays, with a length 28 and with values between –1 and 1, for 3500 evaluations. Hyperopt was run similarly, generating the 56 values between –1 and 1 using the 'hp.uniform' function, and set to evaluate 3500 times using the Tree-structured Parzen Estimator. Paddy was run in both Generational and Population mode, with the Gaussian type set to scaled and with limits of –1 and 1, and randomly propagated in range of the limits with 0.05 in resolution. PFARunner parameters were set where: the random seed number as 250, yt as 15, Qmax as 25, r as 5, and iterations being 30.

Trials run using the multi-feature custom metric were done using the same parameters for those using solely Tversky Similarity. Our custom metric was developed using various methods found in the RDKit library by expanding upon, and modifying terms, the target chemical property function described in the JTVAE manuscript.[86] The target chemical property function was defined by Jaakkola et all as: the octanol-water partition coefficient (*LogP*) of a molecule, minus the molecules Synthetic Accessibility (*SA*) Score and number of cycles with an atom count greater than six (*cycle*) (**eq 7**). For our custom metric, we incorporated both *SA* and *cycle* while introducing: Tversky Similarity (*TS*), fingerprint density (*FD*), rotatable bonds, the number of cycles, and number of on bits (**eq 8**).

Tversky Similarity was calculated using the same parameters as described previously. Fingerprint density was calculated using the RDKit 'FpDensityMorgan3' , which generates Morgan Fingerprints as undefined integer sparse bit vectors with a bit radius of three and returns the quotient of on bits by the number of non-hydrogen atoms in the molecule. Fingerprint density was used with the intention of promoting structurally diverse molecules. Rotatable bonds were enumerated for molecules via RDKit, and used in a conditional manner to calculate a Rotatable Bond Score ($RBS$, see **eq 9**), and penalize long chain and flat molecules. The cycle calculation was expanded as to provide a conditional Cycle Count Score ($CCS$, see **eq 10**) to promote the generation of molecules containing between two and five rings. The RDKit method 'GetMorganFingerprintAsBitVect', which generates explicit bit vectors, was employed with the parameters previously described, with the number of on bits ($mb$) used to penalize molecules with less than 45 on positions as the Bit On Score ($BOS$, see **eq 11**). These individual scores are then used to formulate our custom metric (**eq 8**).

$$f(m) = LogP - SA - cycle \quad \ldots \ldots \ldots \ldots \ldots \ldots \ldots \ldots \ldots (7)$$

$$f^*(m) = TV \cdot FD^2 \cdot BOS \cdot 0.1^{RBS \cdot CCS} \cdot \left(\frac{1}{SA} + cycle\right) \quad \ldots \ldots \ldots \ldots \ldots \ldots \ldots \ldots \ldots (8)$$

$$RBS(m) = \begin{cases} 2 - mr, & \text{if } mr \leq 2 \\ mr - 5, & \text{if } mr \geq 7 \\ 0, & \text{if } 2 < mr < 7 \end{cases} \quad \ldots \ldots \ldots \ldots \ldots \ldots \ldots \ldots \ldots (9)$$

$$CCS(m) = \begin{cases} |mc - 2|, & \text{if } mc \leq 2 \\ |mc - 5|, & \text{if } mc > 5 \\ 0, & \text{if } 2 < mc < 6 \end{cases} \quad \ldots \ldots \ldots \ldots \ldots \ldots \ldots \ldots \ldots (10)$$

$$BOS(m) = \begin{cases} 0.6(mb - 45), & \text{if } mb - 45 < 0 \\ 1, & \text{if } mb \geq 45 \end{cases} \quad \ldots \ldots \ldots \ldots \ldots \ldots \ldots \ldots \ldots (11)$$

$mr$, rotatable bonds in a molecule $m$

$mc$, number of cycles in a molecule $m$

$mb$, number of on bits in an explicit bit vector for a molecule $m$

*Visualization of Latent Space*

Libraries used for analysis and visualization, were UMAP and Matplotlib respectively.[87] The 56 values in the latent vectors generated by Paddy were used as input features. The 56 length vectors were then reduced via UMAP to three components, while setting the number of neighbors to 15 and minimum distance between projected points to 0.5. All other parameters were set to defaults.

## Results and Discussion

### Paddy identifies correct Global Maxima of a Three-Dimensional Bimodal Distribution

We employed a three-dimensional bimodal distribution function with two parameters $(x, y)$ to assess the performance of Paddy to identify the global maximum, out of two maxima (**Figure 2a**). The slope of the global maximum is steeper than that of the local maximum that presents a challenge for global optimization as there is a greater probability that initial sampling will occur near the local maximum with rare events at the global maxima. Both Paddy and Hyperopt were evaluated 100 times, with different starting conditions, with Paddy finding the global maxima 74 times as compared to 13 times for Hyperopt (**Figure 2b-c**). Furthermore, Paddy was 3.4 times faster, on average, than the Bayesian method in Hyperopt. These result support that Paddy sampling may be useful to identify rare events without incurring a high computational cost.

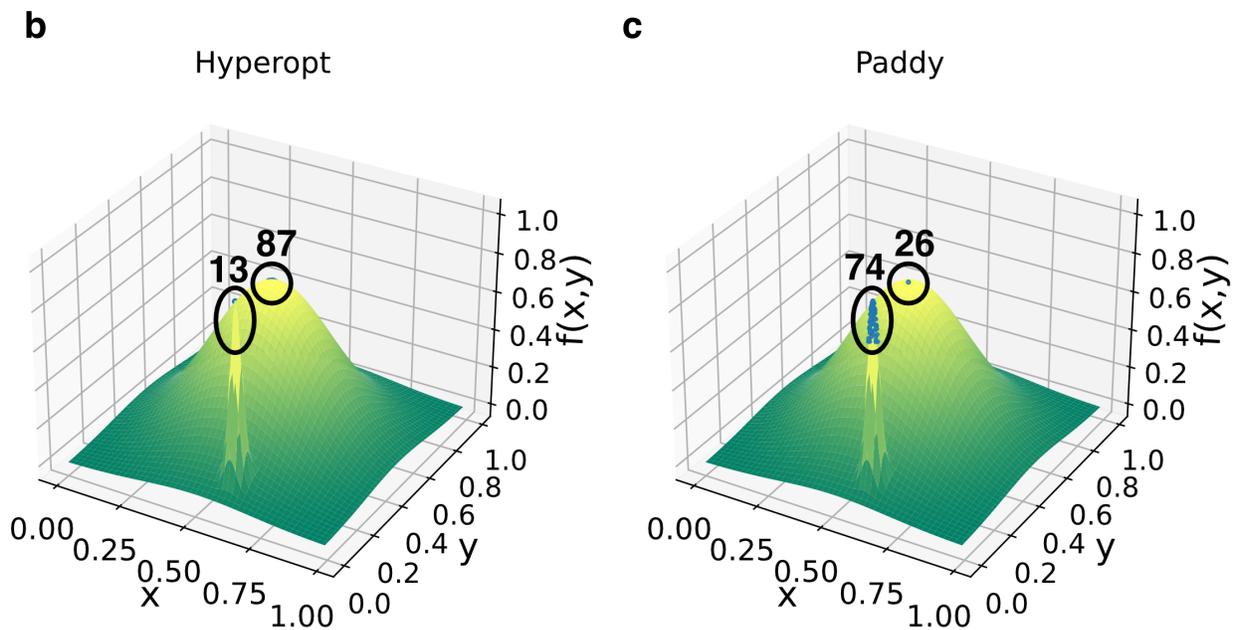

**Figure 2.** Global optimization of a (a) 3D bimodal distribution $f(x, y)$, with a local maxima at (0.5,0.5) and a global maxima at (0.6,0.1). Over 100 runs, Hyperopt (b) successfully identifies the global maxima 13 times, while Paddy (c) identifies the global maxima 74 times respectively without knowledge of the underlying mathematic function.

### Interpolation of the Gramacy & Lee Function using Paddy

To showcase the use of optimization problems to efficiently sample several parameters, we used interpolation of the Gramacy & Lee function[88] using a 65[th] degree trigonometric polynomial, as an example to showcase a possible future application in parameter selection for design of experiments. The performance was evaluated as the mean squared error (MSE) between the $y$ values generated by the 65 fitted polynomial coefficients and the Gramacy & Lee objective function, where $x \in [-0.5, 2.5]$ and with a resolution of 0.001. In order to assess the robustness of the algorithms, we evaluated the performance of both Paddy and Hyperopt optimization for 100 different runs. Paddy displayed both superior performance with lower MSE and runtime when benchmarked against Hyperopt and a random sampling algorithm (**Figure 2a, Table 1**). The MSE performance of Paddy was 0.79 less than Hyperopt, averaging 3.04 and 3.83 respectively. In comparison, the random search algorithm resulted in MSE of 8.69. The best resulting trials for Paddy, Hyperopt and random search are show in comparison to the Gramacy & Lee function (**Figure 2b-d**) and the respective MSE of best fit generated by Paddy was 1.44 lower than Hyperopt (1.94 vs 2.78, **Table 1**). In addition, Paddy was substantially faster with an average runtime of 276 ±13.82 seconds, which is approximately two times faster compared to Hyperopt at 627 ±23.06 seconds and approximately three times faster than random search at 839 ±3.34. These results show that Paddy is an efficient algorithm even when several parameters are sampled for optimization; we expect this algorithm to perform well in design of experiments tasks.

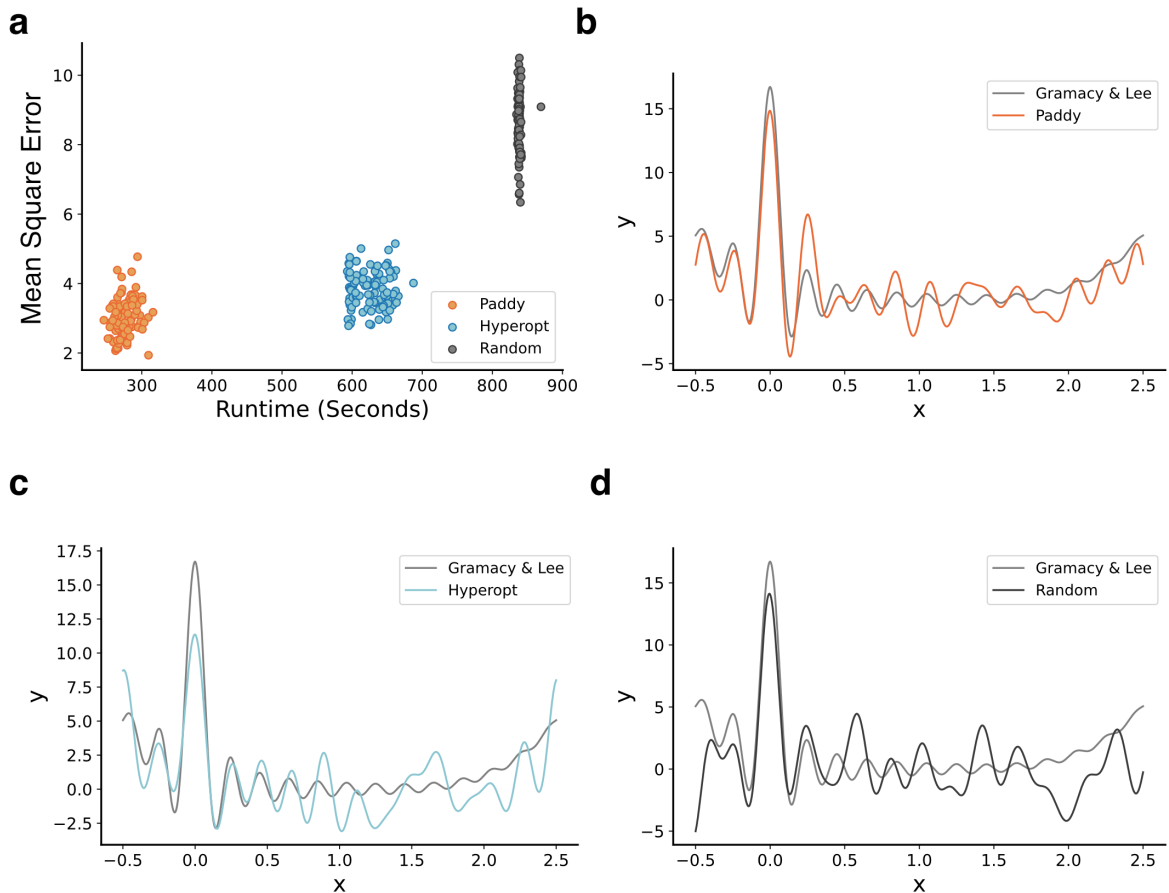

**Figure 3.** Interpolation of the Gramacy & Lee equation using a 65[th] degree trigonometric polynomial.

Scatter plot (a) of mean square error against runtime for Paddy, Hyperopt, and random search, each run 100 times. Plots of the best fitted polynomials for Paddy (b), Hyperopt (c), and random search (d) alongside Gramacy & Lee, with mean squared errors of 1.94, 2.78, and 6.34 respectively.

**Table 1. Gramacy & Lee Interpolation Results**

| Algorithm | Hyperopt | Paddy | Random |
|---|---|---|---|
| best fit (MSE) | 2.78 | 1.94 | 6.34 |
| worst fit (MSE) | 5.15 | 4.77 | 10.50 |
| average fit (MSE) | 3.83 ±0.53 | 3.04 ±0.53 | 8.69 ±0.84 |
| average runtime (s) | 627 ±23.06 | 276 ±13.82 | 839 ±3.34 |

± root mean squared error

### Paddy supports Hyperparameter Optimization of a Multilayer Perceptron

The AI/ML architectures, such as artificial neural networks have been used extensively in cheminformatics, bioinformatics, and computational chemistry/biology in recent years[36,37,48,89]. However, training large AI/ML models present a major challenge for lowering computational costs for training/validation to efficiently select hyperparameters[90] while maintaining performance of the models. To showcase an example of Paddy for efficient use of hyperparameter optimization, we used a multiplayer perceptron (MLP) with two hidden layers (**Figure 4a**). This MLP was designed as a multiclass classifier trained to classify reactions by selecting suitable solvent, such that the reaction inputs were represented as Morgan Fingerprints were trained with the output for one of 32 solvent labels. The average F1 score resulting from 3-fold cross validation was used as the objective function to sample hyperparameters. Specifically, we assessed the performance of Paddy, Hyperopt and random search to select the number of neurons and the dropout rate for the two hidden layers. The number of neurons is an integer in the range of 300-3000 neurons for the first layer and 32-2000 neurons for the second layer. Dropout is a real number between 0 and 1. The ability for Paddy to confine the values used during random initiation was employed to apply these constrains, with dropout values from 0 – 0.5 and lengths of 500-1000 and 32-500 neurons for the first and second hidden layers respectively. To showcase robustness of the method, 100 trials of each method were done and we found both Paddy and Hyperopt outperformed the random sampling algorithm (**Table 2**). The architectures generated by Paddy displayed marginally greater F1 scores than Hyperopt, however, the runtime for Paddy was significantly less compared to hyperopt (**Figure 4b**). With an average runtime of 639 seconds, compared to the 1058 second average runtime for Hyperopt, Paddy was able to optimize the MLP hyperparameters ~1.7 times faster. The ability to optimize hyperparameters in a facile manner is of great importance, and this example showcase the ability for Paddy to be a suitable platform for hyperparameter selection of neural network architectures.

**Table 2. Hyperparameter Optimization Results**

| Algorithm | Hyperopt | Paddy | Random |
|---|---|---|---|
| best F1 | 0.663 | 0.665 | 0.637 |
| worst F1 | 0.640 | 0.645 | 0.588 |
| average F1 | 0.653 ±3.40x10$^{-3}$ | 0.656 ± 3.19x10$^{-3}$ | 0.610 ± 9.52x10$^{-3}$ |
| average runtime (s) | 1058 ±10.71 | 639 ±57.99 | 946 ±19.45 |

± root mean squared error

**a**

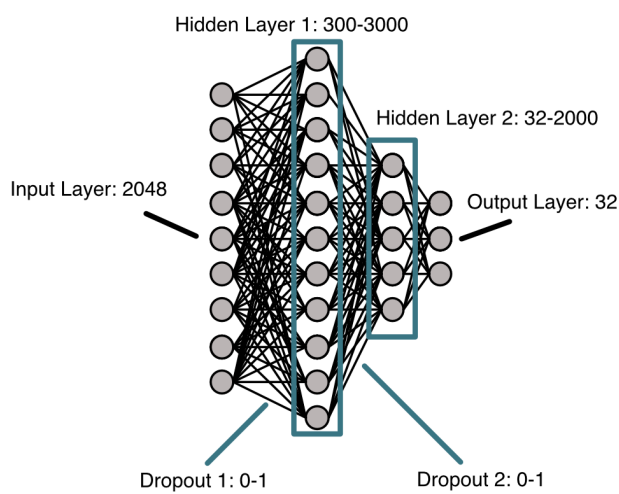

**b**

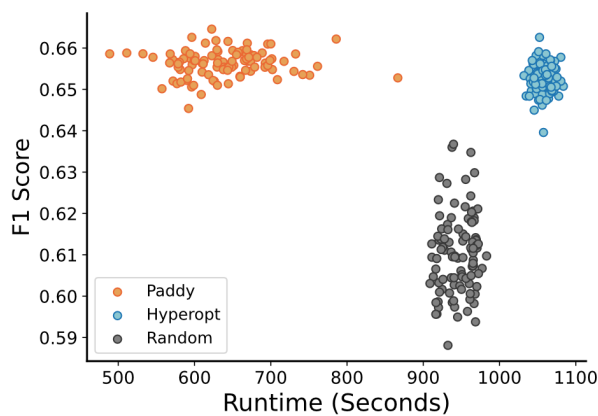

**Figure 4**. Optimization of hyperparameters for a neural network trained to predict solvent for the Morgan fingerprint of reaction components. The architecture of the neural network (a) contained two hidden layers, with the length and dropout of each layer being the objective function variables to optimize. Scatter plot (b) of F1 score against runtime for Paddy, Hyperopt and random search, plotting the greatest

F1 score, as the average from 3-fold cross validation post training, generated over 100 independent runs for each algorithm.

### Sampling Latent Space with Paddy for Targeted Molecule Generation

Another popular application of AI/ML models in chemical sciences is the use of generative neural networks[91], where the model learns the mapping between input and output from random inputs. These models are then used to generate desired outputs satisfying specific conditions such as experimental conditions[92] or molecular structures[45,56,60] based on the mapping of random input to the desired output from the training set. For the task of molecule generation, a popular neural network architecture employed are encoders and decoders[59,60,93,94]. An example of encoder/decoder architecture is an autoencoder, a neural network that is trained to reduce dimensionality of an input and subsequently generate an output, that is ideally, identical to the initial input[95]. The portion of the network tasked with dimension reduction is the encoder network, and the network that reconstructs the input being the decoder network. Transient values feed forward between the encoder and decoder, are often referred to as either latent representations or latent variables. Once trained, autoencoders can then be used in a generative manner by providing a latent representation as input to the decoder network. Furthermore, targeted generation can be conducted via sequential optimization of a latent vector, with examples of this in drug discovery[45,59].

While latent representations have been used for non-generative tasks, the variational autoencoder (VAE) has emerged as an architecture particularly well suited for generative tasks. This is due to the latent variables of VAEs being regularized, as VAEs are trained to optimize the parameters for set of normal distributions and subsequently decode from latent vectors propagated from these distributions. Regularization of the latent space is further reinforced such that the learned distributions are trained to fit a standard normal distribution, with a mean of 0 and standard deviation of 1, in parallel to input reconstruction fidelity. The regularization of the latent distributions results in continuous latent spaces with minimized sparsity. Due to these features, VAEs are well suited for generative tasks, partially as latent space sampling often generates outputs similar in nature to those of neighboring latent features.[96]

To showcase the ability of Paddy to optimally sample latent space vectors, with the goal of target molecule generation, we selected junction tree VAE (JT-VAE)[86]. The JT-VAE functions as a VAE while encoding and decoding molecular graphs with a high degree of reconstructive accuracy. For this case of targeted molecule generation, we utilized Tversky Similarity[97] and our own multi-feature objective function to provide fitness metrics for Paddy and Hyperopt (**Figure 5a**). For the following trials, we used Pazopanib as the target molecule of interest, and a JT-VAE model trained with the ZINC dataset directly taken from the JT-VAE repo.

We used Tversky (Index) Similarity to compare the associated Morgan Fingerprints[83] of generated molecules against Pazopanib, and benchmarked both Paddy and Hyperopt against a random sampling algorithm. Tversky Similarity is the generalized form of Tanimoto Similarity, which are both similarity measurers used to compare sets. Tversky Similarity and Tanimoto Similarity differ where α and β coefficients are set to equal 1 for Tanimoto Similarity and are arbitrary for Tversky Similarity (**eq. 12**). For this instance, the sets are the Morgan Fingerprints, with the bit values for hashed subgraphs being the elements. Set $X$ was the sampled fingerprint while $Y$ was the fingerprint of Pazopanib. The coefficients α and β were set to 0.5 and 0.01 respectively. The low value β was assigned to reduce the penalty for Pazopanib subgraphs not being present in the generated fingerprint.

$$S(X,Y) = \frac{|X \cap Y|}{|X \cap Y| + \alpha|X \setminus Y| + \beta|Y \setminus X|} \quad \ldots \ldots \ldots \ldots \ldots \ldots \ldots \ldots \ldots \ldots (12)$$

The molecule generated by the random sampling algorithm with the greatest fitness was then used as a baseline for comparing the diversity of high similarity molecules generated in turn by Paddy and Hyperopt. This approach for comparing algorithm performance was also employed with the use of our multi-feature objective. To provide further emphasis on drug likeness for generated molecules our custom metric considered in addition to Tversky Similarity: rotatable bonds, the number of cycles, size of cycles, synthetic accessibility, the number of on bits in Morgan Fingerprints, and the number of non-hydrogen atoms (see methods).

Results from generative sampling of latent space using the two metrics described prior indicated that Paddy is well suited for such a task. Paddy generated molecules with greater maximal fitness, less runtime, and a larger population of molecules outperforming the random search solution than Hyperopt. We found Hyperopt quickly optimizes latent space sampling, though plateauing in performance, whereas Paddy avoids early convergence (**Figure 5b-c**). The top scoring molecules generated by both algorithms managed to capture the m-toluenesulfonyl moiety of Pazopanib, however the Tversky Similarity metric rewarded generation of molecules with little chemical diversity as to minimize dissimilarity (**Figure 5d**), which was mitigated by our custom metric (**Figure 5e**). Analysis of the SMILES strings generated by Hyperopt indicates that the algorithm repeatedly samples latent space in the same location after finding a local solution. The convergent behavior of Hyperopt is illustrated by having generated the same solution 249 times using Tversky Similarity (**Table 3**) and 586 times with our custom metric (**Table 4**). The SMILES strings of solutions generated by Paddy and Hyperopt can be found in the supporting information (**Table S1-6**).

Comparing the performance of Paddy, when run in Generational mode versus Population mode, we found the two modes generate differing results while both outperforming Hyperopt. When using Tversky Similarity, incidentally, the two Paddy Modes generated the same number of unique solutions, though Generational mode resulted in a solution of lower similarity with slightly more evaluations and runtime (**Table 3**). Using our custom metric, Generational mode again produced a top solution with a lower score compared to Population mode and with a greater runtime. However, Generational mode yielded nearly twice the number of unique solutions (**Table 4**).

For the optimizations using Tversky Similarity the behavior of the two paddy modes were more so analogous (**Figure 6a**). While the average performance per iteration for both modes was nearly identical, the two diverged in terms of top seed performance. Though generation mode produced solutions sooner than Population mode, Population mode overtook the performance of Generational mode halfway through the run. A greater discrepancy in general behavior was observed between the two modes when using the multi-feature custom metric (**Figure 6b**). Solutions produced by Generational mode displayed a greater average fitness per iteration, and the Generational mode run was only to be bested by Population mode much closer to the end of the run. As the custom-metric accounts for multiple molecular features, this difference in performance may be a result of Population mode being better suited for rapid optimization of relatively smooth response surfaces. Generation mode, however, is inherently more explorative, as it does not sow using the full population of seeds generated during a run. This would lend to the notion of Generational mode being better suited for avoiding repeated sampling of local solutions.

Comparing the generation of solutions by Population and Generational mode, previous insight regarding task specific behavior can be further reinforced. Using Tversky Similarity, the two Paddy modes display analogous behavior, as described in prior, with both optimizing similarity between throughout the run (**Figure 6c**). Both modes display the same general trends in optimization, generating solutions with increasing fitness while followed by discovery of lower scoring solutions. It is interesting however, to note that there were only two identical solution molecules generated by both algorithms. This low frequency of overlap using Tversky Similarity is contrasted by results from using the multi-feature objective function, where various solutions are both identical and, in some cases, generated by both modes during the same Paddy iteration (**Figure 6d**). The overlap in generated solutions would indicate that both Paddy modes sampled latent space in close spatial proximity in part, though with Generational mode having sampled both over a larger area and generated a greater number of solutions. A uniform manifold approximation and projection (UMAP)[87] plot (**Figure 7**) supports this, with Population mode and Generational mode diverging in latent space and Generational mode covering a wider area (Supporting GIF).

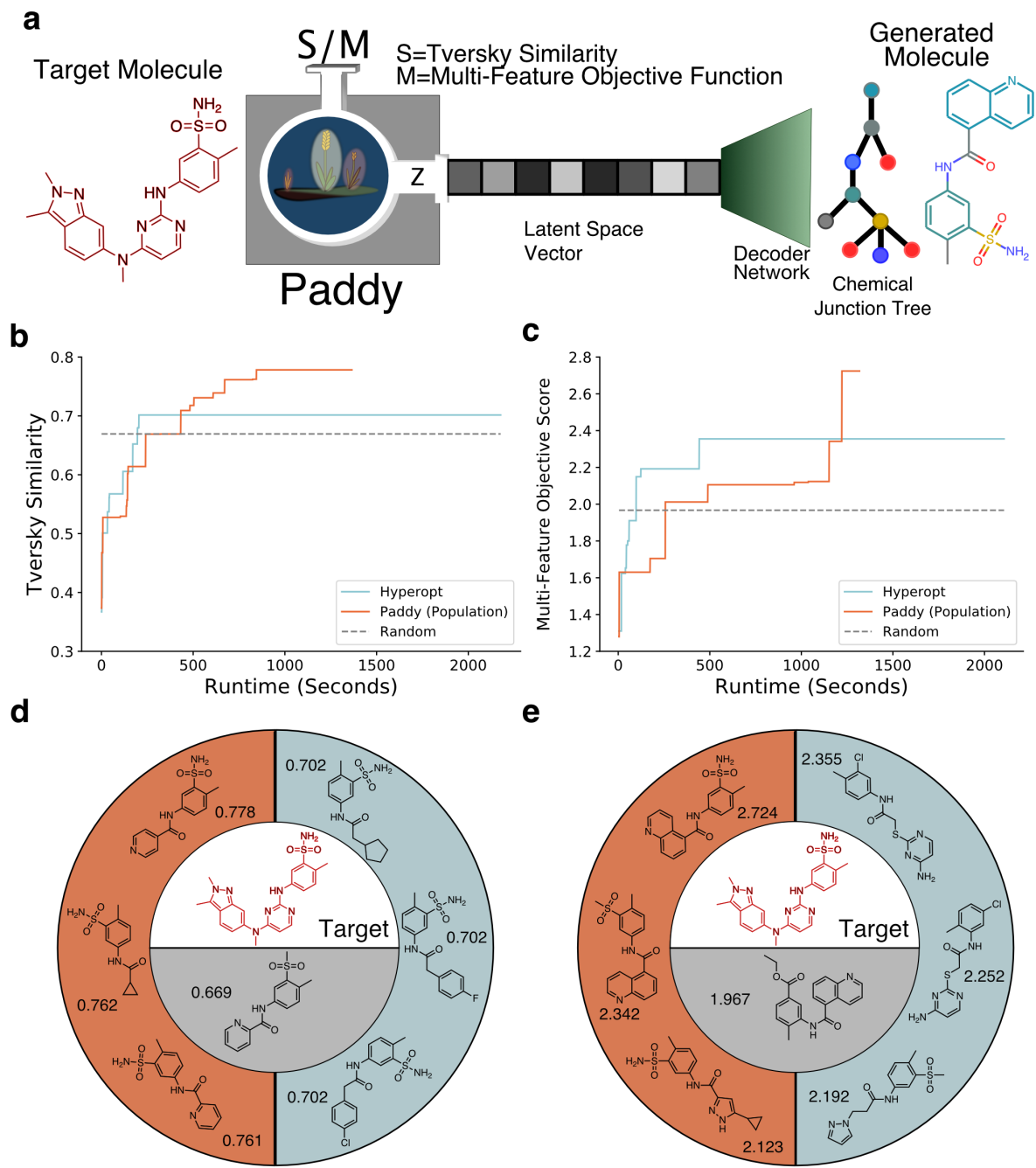

**Figure 5.** Overview of Paddy-JTVAE pipeline (a), where latent space vectors are optimized, such that when decoded to chemical junction trees, generating molecules to maximize an objective function that incorporates a target molecule (Pazopanib). Tversky Similarity (b) and a multifeatured objective function (c) trials are plotted as the running solution over runtime for Paddy and Hyperopt, where the highest scoring random search solution is plotted as a dashed line for comparison. The top three molecules generated by Paddy (orange) and Hyperopt (blue) for the Tversky Similarity (d) and multifeatured objective function (e) trials are displayed with their respective scores, with Pazopanib in red and the random search solution bellow (gray).

**Table 3. Performance using Tversky Similarity as Objective Function**

| Algorithm | Paddy (Population) | Paddy (Generational) | Hyperopt | Random |
|---|---|---|---|---|
| best solution | 0.778 | 0.776 | 0.702 | 0.699 |
| runtime (seconds) | 1365 | 1441 | 2232 | 1171 |
| total evaluations | 4107 | 3571 | 3500 | 3500 |
| unique solutions | 20 | 25 | 14 | ---------- |

**Table 4. Performance When using Custom Multi Feature Objective Function**

| Algorithm | Paddy (Population) | Paddy (Generational) | Hyperopt | Random |
|---|---|---|---|---|
| best solution | 2.724 | 2.265 | 2.355 | 1.967 |
| runtime (seconds) | 1317 | 1849 | 2120 | 1170 |
| total evaluations | 3643 | 5035 | 3500 | 3500 |
| unique solutions | 18 | 33 | 7 | ---------- |

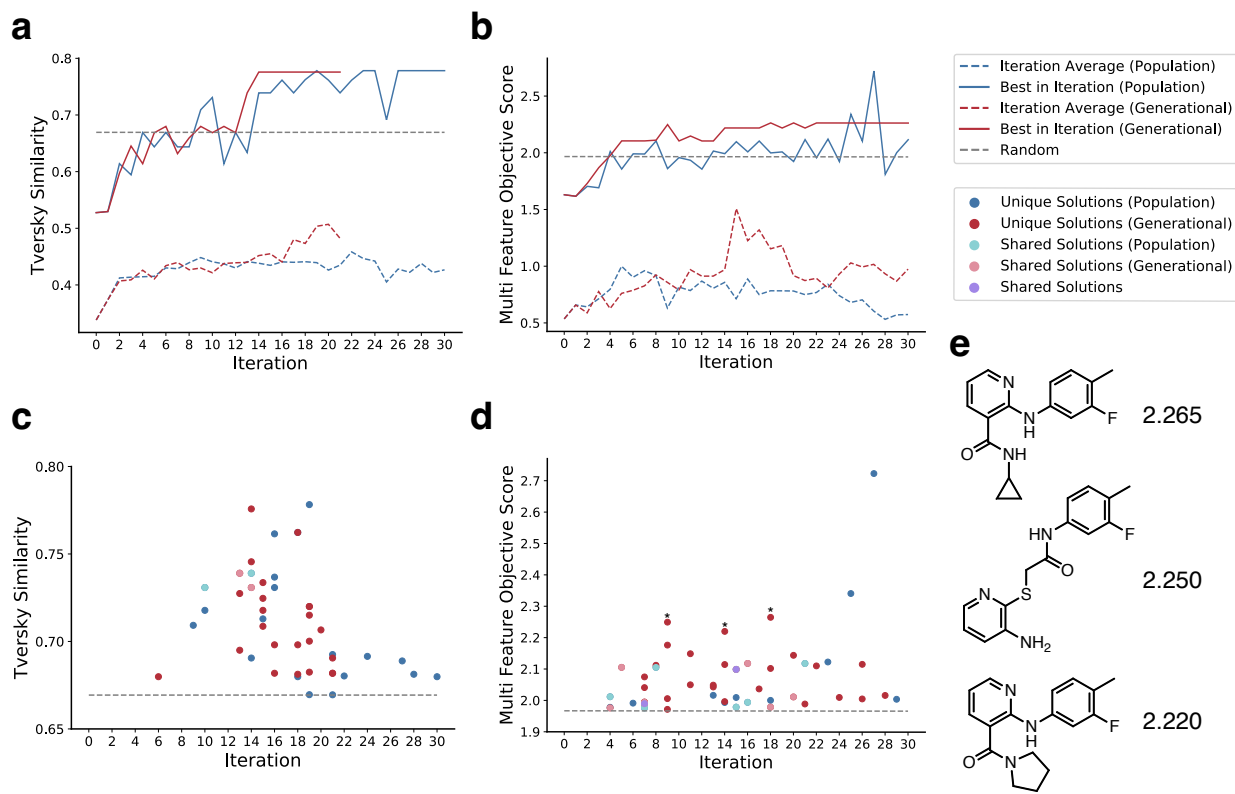

**Figure 6. Comparison of Paddy-JTVAE using population versus generational modes.** Line plots when using Tversky Similarity (a) and the multi-featured objective function (b) depict the performance as the top score evaluated (solid colored) per iteration, and the average performance per (dashed colored) iteration. Scatter plots for the trials using Tversky Similarity (c) and the multi-featured objective function (d) depict the first instance of generating a molecule of a greater score than the maximal performance from random search. Random search performance values for respective metrics are presented as dashed grey lines. The top three molecules generated in generational mode using the multi-feature objective function (e) are displayed with their respective score and denoted with asterisk on the scatter plot.

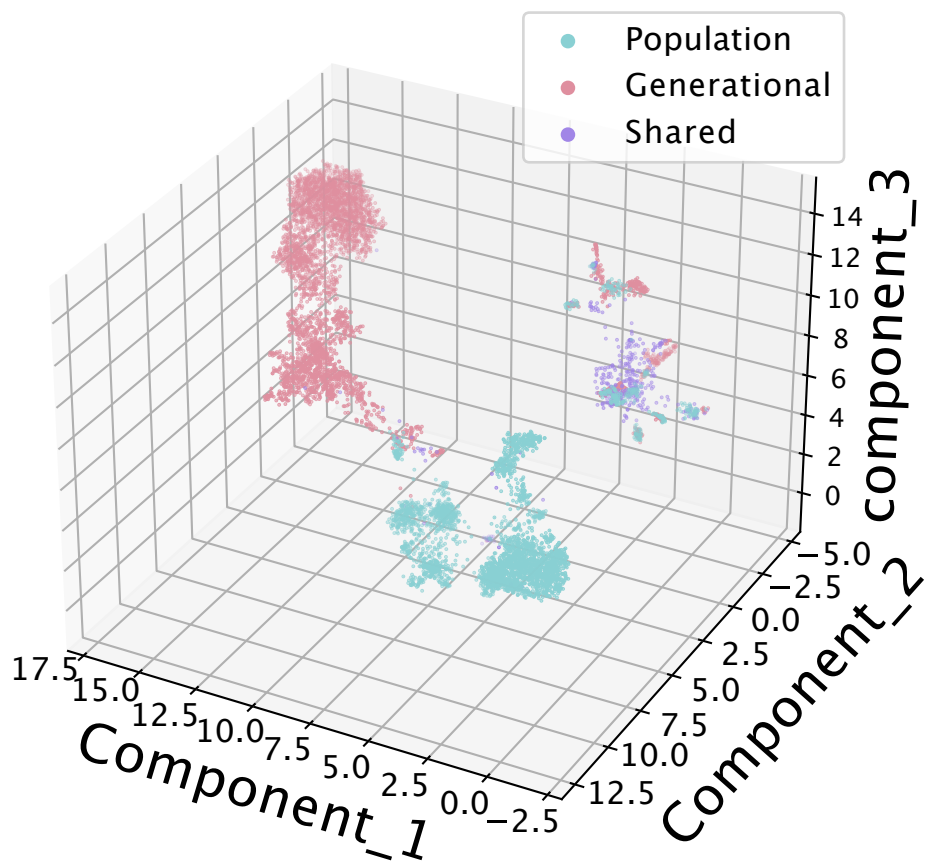

**Figure 7.** Visualization of the latent space sampled by Paddy in both Generational and Population mode when using the multi-feature objective function via UMAP. Shared seeds, due to both modes being initiated with the same random variables, are depicted in purple, where Generational and Population modes are magenta and blue respectively. The greater area sampled by Generational mode and its overlap with Population mode indicates a greater propensity for explorative behavior during optimization, where Population mode appears more exploitive by comparison. This can be visualized best as a rotating gif image provided on GitHub,
https://github.com/chopralab/Paddy_Manuscript_Repo/blob/main/Plotting/JTVAE_Sampling/umap.gif.

## Summary

We introduce an evolutionary algorithm, Paddy, as a python library containing various methods based on the PFA, for facile optimization of numeric parameters for several applications. By considering the spatial distance of parameters and their evaluated performance, Paddy can efficiently optimize a variety of systems without inferring the underlying function. We have benchmarked Paddy against the Tree-structured Parzen Estimator implemented in Hyperopt and we have found Paddy to optimize with less runtime while also avoiding early convergence on a local minimum/maximum. In the context of chemioinformatics, we have shown Paddy to perform well with the tasks of hyperparameter optimization and targeted molecule generation. Additionally, we have investigated the differences in behavior between the native PFA Population mode and our variant, Generational mode, and have shown our variant to be better suited for exploitative optimization while still retaining general performance. We believe that these qualities make Paddy well suited for the optimization of chemical systems of high dimensionality and suitable for tasks such as autonomous closed-loop experimentation and inverse design of drug candidates.

## Conflicts of Interest

Gaurav Chopra is the Director of Merck-Purdue Center for Measurement Science funded by Merck Sharp & Dohme LLC., a subsidiary of Merck & Co., Inc., Rahway, NJ, U.S.A. and the co-founder of Meditati Inc. All other authors declare no competing financial interests.

## Acknowledgements

This work was supported, in part, by the NSF I/UCRC Center for Bioanalytical Metrology (Award 1916991), Purdue Integrative Data Science Institute award, the National Institutes of Health (NIH) award, R01MH128866 by National Institute of Mental Health, and NIH National Center for Advancing Translational Sciences U18TR004146 and ASPIRE Challenge and Reduction-to-Practice awards to G.C. The Purdue University Center for Cancer Research funded by NIH grant P30 CA023168 is also acknowledged. The content is solely the responsibility of the authors and does not necessarily represent the official views of the National Institutes of Health.

**References**

(1) David B. Berkowitz; Mohua Bose; Sungjo Choi. In Situ Screening to Optimize Variables in Organic Reactions, September 6, 2002. https://patents.google.com/patent/US6974665 (accessed 2019-04-15).

(2) Fan, J.; Yi, C.; Lan, X.; Yang, B. Optimization of Synthetic Strategy of 4'4"(5")-Di- *Tert* -Butyldibenzo-18-Crown-6 Using Response Surface Methodology. *Org Process Res Dev* **2013**, *17* (3), 368–374. https://doi.org/10.1021/op3003163.

(3) Sans, V.; Porwol, L.; Dragone, V.; Cronin, L. A Self Optimizing Synthetic Organic Reactor System Using Real-Time in-Line NMR Spectroscopy. *Chem Sci* **2015**, *6* (2), 1258–1264. https://doi.org/10.1039/C4SC03075C.

(4) Welch, C. J.; Regalado, E. L. Estimating Optimal Time for Fast Chromatographic Separations. *J Sep Sci* **2014**, *37* (18), 2552–2558. https://doi.org/10.1002/jssc.201400508.

(5) Glajch, J. L.; Kirkland, J. J.; Squire, K. M.; Minor, J. M. Optimization of Solvent Strength and Selectivity for Reversed-Phase Liquid Chromatography Using an Interactive Mixture-Design

Statistical Technique. *J Chromatogr A* **1980**, *199*, 57–79. https://doi.org/10.1016/S0021-9673(01)91361-5.

(6) Siouffi, A. M.; Phan-Tan-Luu, R. Optimization Methods in Chromatography and Capillary Electrophoresis. *J Chromatogr A* **2000**, *892* (1–2), 75–106. https://doi.org/10.1016/s0021-9673(00)00247-8.

(7) Minenkov, Y.; Singstad, Å.; Occhipinti, G.; Jensen, V. R. The Accuracy of DFT-Optimized Geometries of Functional Transition Metal Compounds: A Validation Study of Catalysts for Olefin Metathesis and Other Reactions in the Homogeneous Phase. *Dalton Transactions* **2012**, *41* (18), 5526. https://doi.org/10.1039/c2dt12232d.

(8) Chen, R.; Christiansen, M. G.; Anikeeva, P. Maximizing Hysteretic Losses in Magnetic Ferrite Nanoparticles via Model-Driven Synthesis and Materials Optimization. *ACS Nano* **2013**, *7* (10), 8990–9000. https://doi.org/10.1021/nn4035266.

(9) Ziaee, A.; Chovan, D.; Lusi, M.; Perry, J. J.; Zaworotko, M. J.; Tofail, S. A. M. Theoretical Optimization of Pore Size and Chemistry in SIFSIX-3-M Hybrid Ultramicroporous Materials. *Cryst Growth Des* **2016**, *16* (7), 3890–3897. https://doi.org/10.1021/acs.cgd.6b00453.

(10) Chen, C.-T.; Chrzan, D. C.; Gu, G. X. Nano-Topology Optimization for Materials Design with Atom-by-Atom Control. *Nat Commun* **2020**, *11* (1), 3745. https://doi.org/10.1038/s41467-020-17570-1.

(11) Fine, J.; Wijewardhane, P. R.; Mohideen, S. D. B.; Smith, K.; Bothe, J. R.; Krishnamachari, Y.; Andrews, A.; Liu, Y.; Chopra, G. Learning Relationships Between Chemical and Physical Stability for Peptide Drug Development. *Pharm Res* **2023**, *40* (3), 701–710. https://doi.org/10.1007/s11095-023-03475-3.

(12) Bushnell, E. A. C.; Burns, T. D.; Boyd, R. J. The One-Electron Oxidation of a Dithiolate Molecule: The Importance of Chemical Intuition. *J Chem Phys* **2014**, *140* (18), 18A519. https://doi.org/10.1063/1.4867537.

(13) Keserű, G. M.; Soós, T.; Kappe, C. O. Anthropogenic Reaction Parameters – the Missing Link between Chemical Intuition and the Available Chemical Space. *Chem. Soc. Rev.* **2014**, *43* (15), 5387–5399. https://doi.org/10.1039/C3CS60423C.

(14) Lombardino, J. G.; Lowe, J. A. The Role of the Medicinal Chemist in Drug Discovery--Then and Now. *Nat Rev Drug Discov* **2004**, *3* (10), 853–862. https://doi.org/10.1038/nrd1523.

(15) Yao, K.; Liu, M.; Zheng, Z.; Shih, T.; Xie, J.; Sun, H.; Chen, Z. Automatic Shimming Method Using Compensation of Magnetic Susceptibilities and Adaptive Simplex for Low-Field NMR. *IEEE Trans Instrum Meas* **2021**, *70*, 1–12. https://doi.org/10.1109/TIM.2021.3074951.

(16) Zhang, J.; Gonzalez, E.; Hestilow, T.; Haskins, W.; Huang, Y. Review of Peak Detection Algorithms in Liquid-Chromatography-Mass Spectrometry. *Curr Genomics* **2009**, *10* (6), 388–401. https://doi.org/10.2174/138920209789177638.

(17) Evans, B. R.; Yeung, L.; Beck, A. G.; Li, A.; Lee, D. H.; Chopra, G.; Bateman, K. P. Automated Bioanalytical Workflow for Ligand Binding Based Pharmacokinetic Assay Development. *ArXiv* **2022**. https://doi.org/10.26434/chemrxiv-2022-b3gn1.

(18) Wu, Q.; Yang, W. Algebraic Equation and Iterative Optimization for the Optimized Effective Potential in Density Functional Theory. *J Theor Comput Chem* **2003**, *02* (04), 627–638. https://doi.org/10.1142/S0219633603000690.

(19) Kong, J.; Eason, J. P.; Chen, X.; Biegler, L. T. Operational Optimization of Polymerization Reactors with Computational Fluid Dynamics and Embedded Molecular Weight Distribution Using the Iterative Surrogate Model Method. *Ind Eng Chem Res* **2020**, *59* (19), 9165–9179. https://doi.org/10.1021/acs.iecr.0c00367.

(20) Ostrovsky, G. M.; Ziyatdinov, N. N.; Lapteva, T. V.; Silvestrova, A. Optimization of Chemical Process Design with Chance Constraints by an Iterative Partitioning Approach. *Ind Eng Chem Res* **2015**, *54* (13), 3412–3429. https://doi.org/10.1021/ie5048016.

(21) Varela, R.; Walters, W. P.; Goldman, B. B.; Jain, A. N. Iterative Refinement of a Binding Pocket Model: Active Computational Steering of Lead Optimization. *J Med Chem* **2012**, *55* (20), 8926–8942. https://doi.org/10.1021/jm301210j.

(22) Qian, F.; Sun, F.; Zhong, W.; Luo, N. Dynamic Optimization of Chemical Engineering Problems Using a Control Vector Parameterization Method with an Iterative Genetic Algorithm. *Engineering Optimization* **2013**, *45* (9), 1129–1146. https://doi.org/10.1080/0305215X.2012.720683.

(23) Zheng, J.; Frisch, M. J. Efficient Geometry Minimization and Transition Structure Optimization Using Interpolated Potential Energy Surfaces and Iteratively Updated Hessians. *J Chem Theory Comput* **2017**, *13* (12), 6424–6432. https://doi.org/10.1021/acs.jctc.7b00719.

(24) Zhang, B.; Chen, D.; Zhao, W. Iterative Ant-Colony Algorithm and Its Application to Dynamic Optimization of Chemical Process. *Comput Chem Eng* **2005**, *29* (10), 2078–2086. https://doi.org/10.1016/j.compchemeng.2005.05.020.

(25) Li, D.-W.; Brüschweiler, R. Iterative Optimization of Molecular Mechanics Force Fields from NMR Data of Full-Length Proteins. *J Chem Theory Comput* **2011**, *7* (6), 1773–1782. https://doi.org/10.1021/ct200094b.

(26) Pantazes, R. J.; Grisewood, M. J.; Li, T.; Gifford, N. P.; Maranas, C. D. The Iterative Protein Redesign and Optimization (IPRO) Suite of Programs. *J Comput Chem* **2015**, *36* (4), 251–263. https://doi.org/10.1002/jcc.23796.

(27) Farkas, Ö.; Schlegel, H. B. Methods for Optimizing Large Molecules. Part III. An Improved Algorithm for Geometry Optimization Using Direct Inversion in the Iterative Subspace (GDIIS). *Phys. Chem. Chem. Phys.* **2002**, *4* (1), 11–15. https://doi.org/10.1039/B108658H.

(28) Piris, M.; Ugalde, J. M. Iterative Diagonalization for Orbital Optimization in Natural Orbital Functional Theory. *J Comput Chem* **2009**, *30* (13), 2078–2086. https://doi.org/10.1002/jcc.21225.

(29) Tang, J.; Egiazarian, K.; Golbabaee, M.; Davies, M. The Practicality of Stochastic Optimization in Imaging Inverse Problems. *IEEE Trans Comput Imaging* **2019**, *6*, 1471–1485. https://doi.org/10.1109/TCI.2020.3032101.

(30) Farasat, E.; Huang, B. Deterministic vs. Stochastic Performance Assessment of Iterative Learning Control for Batch Processes. *AIChE Journal* **2013**, *59* (2), 457–464. https://doi.org/10.1002/aic.13840.

(31) Philbrick, C. R.; Kitanidis, P. K. Limitations of Deterministic Optimization Applied to Reservoir Operations. *J Water Resour Plan Manag* **1999**, *125* (3), 135–142. https://doi.org/10.1061/(ASCE)0733-9496(1999)125:3(135).

(32) Pool, M.; Carrera, J.; Alcolea, A.; Bocanegra, E. M. A Comparison of Deterministic and Stochastic Approaches for Regional Scale Inverse Modeling on the Mar Del Plata Aquifer. *J Hydrol (Amst)* **2015**, *531*, 214–229. https://doi.org/10.1016/j.jhydrol.2015.09.064.

(33) Kingma, D. P.; Ba, J. Adam: A Method for Stochastic Optimization. *ArXiv* **2014**.

(34) Zhang, J.; Dolg, M. ABCluster: The Artificial Bee Colony Algorithm for Cluster Global Optimization. *Physical Chemistry Chemical Physics* **2015**, *17* (37), 24173–24181. https://doi.org/10.1039/C5CP04060D.

(35) Alba, E.; Tomassini, M. Parallelism and Evolutionary Algorithms. *IEEE Transactions on Evolutionary Computation* **2002**, *6* (5), 443–462. https://doi.org/10.1109/TEVC.2002.800880.

(36) Goh, G. B.; Hodas, N. O.; Vishnu, A. Deep Learning for Computational Chemistry. *J Comput Chem* **2017**, *38* (16), 1291–1307. https://doi.org/10.1002/jcc.24764.

(37) Mater, A. C.; Coote, M. L. Deep Learning in Chemistry. *J Chem Inf Model* **2019**, *59* (6), 2545–2559. https://doi.org/10.1021/acs.jcim.9b00266.

(38) Liu, B.; Ramsundar, B.; Kawthekar, P.; Shi, J.; Gomes, J.; Luu Nguyen, Q.; Ho, S.; Sloane, J.; Wender, P.; Pande, V. Retrosynthetic Reaction Prediction Using Neural Sequence-to-Sequence Models. *ACS Cent Sci* **2017**, *3* (10), 1103–1113. https://doi.org/10.1021/acscentsci.7b00303.

(39) Zhou, Z.; Li, X.; Zare, R. N. Optimizing Chemical Reactions with Deep Reinforcement Learning. *ACS Cent Sci* **2017**, *3* (12), 1337–1344. https://doi.org/10.1021/acscentsci.7b00492.

(40) Cortés-Borda, D.; Wimmer, E.; Gouilleux, B.; Barré, E.; Oger, N.; Goulamaly, L.; Peault, L.; Charrier, B.; Truchet, C.; Giraudeau, P.; Rodriguez-Zubiri, M.; Le Grognec, E.; Felpin, F.-X. An Autonomous Self-Optimizing Flow Reactor for the Synthesis of Natural Product Carpanone. *J Org Chem* **2018**, *83* (23), 14286–14299. https://doi.org/10.1021/acs.joc.8b01821.

(41) Wei, J. N.; Duvenaud, D.; Aspuru-Guzik, A. Neural Networks for the Prediction of Organic Chemistry Reactions. *ACS Cent Sci* **2016**, *2* (10), 725–732. https://doi.org/10.1021/acscentsci.6b00219.

(42) Gao, H.; Struble, T. J.; Coley, C. W.; Wang, Y.; Green, W. H.; Jensen, K. F. Using Machine Learning To Predict Suitable Conditions for Organic Reactions. *ACS Cent Sci* **2018**, *4* (11), 1465–1476. https://doi.org/10.1021/acscentsci.8b00357.

(43) Zahrt, A. F.; Henle, J. J.; Rose, B. T.; Wang, Y.; Darrow, W. T.; Denmark, S. E. Prediction of Higher-Selectivity Catalysts by Computer-Driven Workflow and Machine Learning. *Science* **2019**, *363* (6424). https://doi.org/10.1126/science.aau5631.

(44) Li, Z.; Wang, S.; Chin, W. S.; Achenie, L. E.; Xin, H. High-Throughput Screening of Bimetallic Catalysts Enabled by Machine Learning. *J Mater Chem A Mater* **2017**, *5* (46), 24131–24138. https://doi.org/10.1039/C7TA01812F.

(45) Putin, E.; Asadulaev, A.; Ivanenkov, Y.; Aladinskiy, V.; Sanchez-Lengeling, B.; Aspuru-Guzik, A.; Zhavoronkov, A. Reinforced Adversarial Neural Computer for *de Novo* Molecular Design. *J Chem Inf Model* **2018**, *58* (6), 1194–1204. https://doi.org/10.1021/acs.jcim.7b00690.

(46) Ying Liu. Drug Design by Machine Learning: Ensemble Learning for QSAR Modeling. In *Fourth International Conference on Machine Learning and Applications (ICMLA'05)*; IEEE, 2005; pp 187–193. https://doi.org/10.1109/ICMLA.2005.25.

(47) Altae-Tran, H.; Ramsundar, B.; Pappu, A. S.; Pande, V. Low Data Drug Discovery with One-Shot Learning. *ACS Cent Sci* **2017**, *3* (4), 283–293. https://doi.org/10.1021/acscentsci.6b00367.

(48) Lo, Y.-C.; Rensi, S. E.; Torng, W.; Altman, R. B. Machine Learning in Chemoinformatics and Drug Discovery. *Drug Discov Today* **2018**, *23* (8), 1538–1546. https://doi.org/10.1016/j.drudis.2018.05.010.

(49) Zhou, C.; Bowler, L. D.; Feng, J. A Machine Learning Approach to Explore the Spectra Intensity Pattern of Peptides Using Tandem Mass Spectrometry Data. *BMC Bioinformatics* **2008**, *9* (1), 325. https://doi.org/10.1186/1471-2105-9-325.

(50) Liu, J.; Zhang, J.; Luo, Y.; Yang, S.; Wang, J.; Fu, Q. Mass Spectral Substance Detections Using Long Short-Term Memory Networks. *IEEE Access* **2019**, *7*, 10734–10744. https://doi.org/10.1109/ACCESS.2019.2891548.

(51) Fine, J. A.; Rajasekar, A. A.; Jethava, K. P.; Chopra, G. Spectral Deep Learning for Prediction and Prospective Validation of Functional Groups. *Chem Sci* **2020**, *11* (18), 4618–4630. https://doi.org/10.1039/C9SC06240H.

(52) Bouwmeester, R.; Martens, L.; Degroeve, S. Comprehensive and Empirical Evaluation of Machine Learning Algorithms for Small Molecule LC Retention Time Prediction. *Anal Chem* **2019**, *91* (5), 3694–3703. https://doi.org/10.1021/acs.analchem.8b05820.

(53) Liu, Y.-B.; Yang, J.-Y.; Xin, G.-M.; Liu, L.-H.; Csányi, G.; Cao, B.-Y. Machine Learning Interatomic Potential Developed for Molecular Simulations on Thermal Properties of β-Ga2O3. *J Chem Phys* **2020**, *153* (14), 144501. https://doi.org/10.1063/5.0027643.

(54) Mittal, S.; Shukla, D. Recruiting Machine Learning Methods for Molecular Simulations of Proteins. *Mol Simul* **2018**, *44* (11), 891–904. https://doi.org/10.1080/08927022.2018.1448976.

(55) Westermayr, J.; Gastegger, M.; Menger, M. F. S. J.; Mai, S.; González, L.; Marquetand, P. Machine Learning Enables Long Time Scale Molecular Photodynamics Simulations. *Chem Sci* **2019**, *10* (35), 8100–8107. https://doi.org/10.1039/C9SC01742A.

(56) Sanchez-Lengeling, B.; Aspuru-Guzik, A. Inverse Molecular Design Using Machine Learning: Generative Models for Matter Engineering. *Science (1979)* **2018**, *361* (6400), 360–365. https://doi.org/10.1126/science.aat2663.

(57) Kim, K.; Kang, S.; Yoo, J.; Kwon, Y.; Nam, Y.; Lee, D.; Kim, I.; Choi, Y.-S.; Jung, Y.; Kim, S.; Son, W.-J.; Son, J.; Lee, H. S.; Kim, S.; Shin, J.; Hwang, S. Deep-Learning-Based Inverse Design Model for Intelligent Discovery of Organic Molecules. *NPJ Comput Mater* **2018**, *4* (1), 67. https://doi.org/10.1038/s41524-018-0128-1.

(58) Benhenda, M. ChemGAN Challenge for Drug Discovery: Can AI Reproduce Natural Chemical Diversity? *ArXiv* **2017**.

(59) Sattarov, B.; Baskin, I. I.; Horvath, D.; Marcou, G.; Bjerrum, E. J.; Varnek, A. De Novo Molecular Design by Combining Deep Autoencoder Recurrent Neural Networks with Generative Topographic Mapping. *J Chem Inf Model* **2019**, *59* (3), 1182–1196. https://doi.org/10.1021/acs.jcim.8b00751.

(60) Kang, S.; Cho, K. Conditional Molecular Design with Deep Generative Models. *J Chem Inf Model* **2019**, *59* (1), 43–52. https://doi.org/10.1021/acs.jcim.8b00263.

(61) Kusne, A. G.; Yu, H.; Wu, C.; Zhang, H.; Hattrick-Simpers, J.; DeCost, B.; Sarker, S.; Oses, C.; Toher, C.; Curtarolo, S.; Davydov, A. V.; Agarwal, R.; Bendersky, L. A.; Li, M.; Mehta, A.; Takeuchi, I. On-the-Fly Closed-Loop Materials Discovery via Bayesian Active Learning. *Nat Commun* **2020**, *11* (1), 5966. https://doi.org/10.1038/s41467-020-19597-w.

(62) Liu, Y.; Yang, J.; Vasudevan, R. K.; Kelley, K. P.; Ziatdinov, M.; Kalinin, S. V.; Ahmadi, M. Exploring the Relationship of Microstructure and Conductivity in Metal Halide Perovskites via Active Learning-Driven Automated Scanning Probe Microscopy. *J Phys Chem Lett* **2023**, *14* (13), 3352–3359. https://doi.org/10.1021/acs.jpclett.3c00223.

(63) Eyke, N. S.; Green, W. H.; Jensen, K. F. Iterative Experimental Design Based on Active Machine Learning Reduces the Experimental Burden Associated with Reaction Screening. *React Chem Eng* **2020**, *5* (10), 1963–1972. https://doi.org/10.1039/D0RE00232A.

(64) Shields, B. J.; Stevens, J.; Li, J.; Parasram, M.; Damani, F.; Alvarado, J. I. M.; Janey, J. M.; Adams, R. P.; Doyle, A. G. Bayesian Reaction Optimization as a Tool for Chemical Synthesis. *Nature* **2021**, *590* (7844), 89–96. https://doi.org/10.1038/s41586-021-03213-y.

(65) Torres, J. A. G.; Lau, S. H.; Anchuri, P.; Stevens, J. M.; Tabora, J. E.; Li, J.; Borovika, A.; Adams, R. P.; Doyle, A. G. A Multi-Objective Active Learning Platform and Web App for Reaction Optimization. *J Am Chem Soc* **2022**, *144* (43), 19999–20007. https://doi.org/10.1021/jacs.2c08592.

(66) Capecchi, A.; Zhang, A.; Reymond, J.-L. Populating Chemical Space with Peptides Using a Genetic Algorithm. *J Chem Inf Model* **2020**, *60* (1), 121–132. https://doi.org/10.1021/acs.jcim.9b01014.

(67) Holmes, N.; Akien, G. R.; Blacker, A. J.; Woodward, R. L.; Meadows, R. E.; Bourne, R. A. Self-Optimisation of the Final Stage in the Synthesis of EGFR Kinase Inhibitor AZD9291 Using an Automated Flow Reactor. *React Chem Eng* **2016**, *1* (4), 366–371. https://doi.org/10.1039/C6RE00059B.

(68) Bertsimas, D.; Tsitsiklis, J. Simulated Annealing. *Statistical Science* **1993**, *8* (1). https://doi.org/10.1214/ss/1177011077.

(69) Katoch, S.; Chauhan, S. S.; Kumar, V. A Review on Genetic Algorithm: Past, Present, and Future. *Multimed Tools Appl* **2021**, *80* (5), 8091–8126. https://doi.org/10.1007/s11042-020-10139-6.

(70) Gendreau, M.; Potvin, J.-Y. Tabu Search. In *Search Methodologies*; Springer US: Boston, MA; pp 165–186. https://doi.org/10.1007/0-387-28356-0_6.

(71) Storey, C. Applications of a Hill Climbing Method of Optimization. *Chem Eng Sci* **1962**, *17* (1), 45–52. https://doi.org/10.1016/0009-2509(62)80005-0.

(72) Kennedy, J.; Eberhart, R. Particle Swarm Optimization. In *Proceedings of ICNN'95 - International Conference on Neural Networks*; IEEE, 1995; Vol. 4, pp 1942–1948. https://doi.org/10.1109/ICNN.1995.488968.

(73) Mockus, J. *Bayesian Approach to Global Optimization*; Mathematics and Its Applications; Springer Netherlands: Dordrecht, 1989; Vol. 37. https://doi.org/10.1007/978-94-009-0909-0.

(74) Ueno, T.; Rhone, T. D.; Hou, Z.; Mizoguchi, T.; Tsuda, K. COMBO: An Efficient Bayesian Optimization Library for Materials Science. *Materials Discovery* **2016**, *4*, 18–21. https://doi.org/10.1016/j.md.2016.04.001.

(75) Griffiths, R.-R.; Hernández-Lobato, J. M. Constrained Bayesian Optimization for Automatic Chemical Design Using Variational Autoencoders. *Chem Sci* **2020**, *11* (2), 577–586. https://doi.org/10.1039/C9SC04026A.

(76) Häse, F.; Roch, L. M.; Kreisbeck, C.; Aspuru-Guzik, A. Phoenics: A Bayesian Optimizer for Chemistry. *ACS Cent Sci* **2018**, *4* (9), 1134–1145. https://doi.org/10.1021/acscentsci.8b00307.

(77) Häse, F.; Aldeghi, M.; Hickman, R. J.; Roch, L. M.; Aspuru-Guzik, A. G<scp>ryffin</Scp> : An Algorithm for Bayesian Optimization of Categorical Variables Informed by Expert Knowledge. *Appl Phys Rev* **2021**, *8* (3), 031406. https://doi.org/10.1063/5.0048164.

(78) Premaratne, U.; Samarabandu, J.; Sidhu, T. A New Biologically Inspired Optimization Algorithm. In *2009 International Conference on Industrial and Information Systems (ICIIS)*; IEEE, 2009; pp 279–284. https://doi.org/10.1109/ICIINFS.2009.5429852.

(79) Bergstra, J.; Bardenet, R.; Bengio, Y.; Kégl, B. Algorithms for Hyper-Parameter Optimization. *Adv Neural Inf Process Syst* **2011**, *24*.

(80) Bergstra, J.; Komer, B.; Eliasmith, C.; Yamins, D.; Cox, D. D. Hyperopt: A Python Library for Model Selection and Hyperparameter Optimization. *Comput Sci Discov* **2015**, *8* (1), 014008. https://doi.org/10.1088/1749-4699/8/1/014008.

(81) Daniel Lowe. Extraction of Chemical Structures and Reactions from the Literature. Thesis, University of Cambridge, 2012. https://doi.org/https://doi.org/10.17863/CAM.16293.

(82) Landrum, G. A. RDKit: Open-Source Cheminformatics. https://www.rdkit.org (accessed 2023-03-08).

(83) Morgan, H. L. The Generation of a Unique Machine Description for Chemical Structures-A Technique Developed at Chemical Abstracts Service. *J Chem Doc* **1965**, *5* (2), 107–113. https://doi.org/10.1021/c160017a018.

(84) Chollet, F. Keras. **2015**.

(85) Pedregosa, F.; Varoquaux, G.; Gramfort, A.; Michel, V.; Thirion, B.; Grisel, O.; Blondel, M.; Müller, A.; Nothman, J.; Louppe, G.; Prettenhofer, P.; Weiss, R.; Dubourg, V.; Vanderplas, J.; Passos, A.; Cournapeau, D.; Brucher, M.; Perrot, M.; Duchesnay, É. Scikit-Learn: Machine Learning in Python. *Journal of Machine Learning Research* **2012**, *12* (85), 2825–2830.

(86) Jin, W.; Barzilay, R.; Jaakkola, T. Junction Tree Variational Autoencoder for Molecular Graph Generation. **2018**.

(87) McInnes, L.; Healy, J.; Melville, J. UMAP: Uniform Manifold Approximation and Projection for Dimension Reduction. *ArXiv* **2018**.

(88) Gramacy, R. B.; Lee, H. K. H. Cases for the Nugget in Modeling Computer Experiments. *Stat Comput* **2012**, *22* (3), 713–722. https://doi.org/10.1007/s11222-010-9224-x.

(89) Tkatchenko, A. Machine Learning for Chemical Discovery. *Nat Commun* **2020**, *11* (1), 4125. https://doi.org/10.1038/s41467-020-17844-8.

(90) Claesen, M.; De Moor, B. Hyperparameter Search in Machine Learning. *ArXiv* **2015**.

(91) White, D.; Wilson, R. C. Generative Models for Chemical Structures. *J Chem Inf Model* **2010**, *50* (7), 1257–1274. https://doi.org/10.1021/ci9004089.

(92) Bort, W.; Baskin, I. I.; Gimadiev, T.; Mukanov, A.; Nugmanov, R.; Sidorov, P.; Marcou, G.; Horvath, D.; Madzhidov, T.; Varnek, A. Discovery of Novel Chemical Reactions by Deep Generative Recurrent Neural Network. *ChemRxiv*. ChemRxiv January 17, 2020. https://doi.org/10.26434/chemrxiv.11635929.v1.

(93) Skalic, M.; Jiménez, J.; Sabbadin, D.; De Fabritiis, G. Shape-Based Generative Modeling for de Novo Drug Design. *J Chem Inf Model* **2019**, *59* (3), 1205–1214. https://doi.org/10.1021/acs.jcim.8b00706.

(94) Born, J.; Manica, M.; Oskooei, A.; Cadow, J.; Markert, G.; Rodríguez Martínez, M. PaccMannRL: De Novo Generation of Hit-like Anticancer Molecules from Transcriptomic Data via Reinforcement Learning. *iScience* **2021**, *24* (4), 102269. https://doi.org/10.1016/j.isci.2021.102269.

(95) Learning, P. B. B. T.-P. of I. W. on U. and T. Autoencoders, Unsupervised Learning, and Deep Architectures. In *Proceedings of the 2011 International Conference on Unsupervised and Transfer Learning Workshop - Volume 27*; Guyon, I., Dror, G., Lemaire, V., Taylor, G., Silver, D., Eds.; PMLR, 2012; pp 37–49.

(96) Zhai, J.; Zhang, S.; Chen, J.; He, Q. Autoencoder and Its Various Variants. In *2018 IEEE International Conference on Systems, Man, and Cybernetics (SMC)*; IEEE, 2018; pp 415–419. https://doi.org/10.1109/SMC.2018.00080.

(97) Kunimoto, R.; Vogt, M.; Bajorath, J. Maximum Common Substructure-Based Tversky Index: An Asymmetric Hybrid Similarity Measure. *J Comput Aided Mol Des* **2016**, *30* (7), 523–531. https://doi.org/10.1007/s10822-016-9935-y.

**Paddy: Evolutionary Optimization Algorithm for Chemical Systems and Spaces**

Armen Beck[1], Jonathan Fine[1], Gaurav Chopra[1,2,3,4,5,6]

[1]Department of Chemistry and Computer Science (*by courtesy*), Purdue University, 720 Clinic Drive, West Lafayette, IN 47907

[2]Purdue Institute for Drug Discovery, West Lafayette, IN 47907

[3]Purdue Center for Cancer Research, West Lafayette, IN 47907

[4]Purdue Institute for Inflammation, Immunology and Infectious Disease, West Lafayette, IN 47907

[5]Purdue Institute for Integrative Neuroscience, West Lafayette, IN 47907

[6]Regenstrief Center for Healthcare Engineering, West Lafayette, IN 47907

*Corresponding author email – gchopra@purdue.edu

# Contents

**Table S1.** Paddy (Population) Solutions using Tversky Similarity as Objective Function

| SMILE string | score | frequency |
|---|---|---|
| Cc1ccc(NC(=O)c2ccncc2)cc1S(N)(=O)=O | 0.778210116732 | 8 |
| Cc1ccc(NC(=O)C2CC2)cc1S(N)(=O)=O | 0.762300762301 | 1 |
| Cc1ccc(NC(=O)c2ccccn2)cc1S(N)(=O)=O | 0.761498629302 | 15 |
| α Cc1ccc(NC(=O)c2cccnc2N)cc1S(N)(=O)=O | 0.738989062962 | 12 |
| Cc1ccc(NC(=O)c2ccc(F)cc2)cc1S(N)(=O)=O | 0.736771600804 | 2 |
| Cc1ccc(NC(=O)c2cccnc2)cc1S(N)(=O)=O | 0.730816077954 | 16 |
| α Cc1ccc(NC(=O)c2ccc(N)cn2)cc1S(N)(=O)=O | 0.730816077954 | 1 |
| Cc1ccc(NC(=O)c2ccc3ccccc3n2)cc1S(N)(=O)=O | 0.717772035601 | 3 |
| Cc1ccc(NC(=O)c2ccccc2F)cc1S(N)(=O)=O | 0.712896953986 | 3 |
| Cc1ccc(NC(=O)c2ccnc(N(C)C)c2)cc1S(C)(=O)=O | 0.709219858156 | 1 |
| Cc1cc(NC(=O)c2ccncc2)ccc1S(N)(=O)=O | 0.692520775623 | 3 |
| Cc1ccc(NC(=O)c2ccncc2)cc1S(=O)(=O)N(C)C | 0.691471847218 | 7 |
| Cc1ccc(NC(=O)c2ccnc(N)c2)cc1S(C)(=O)=O | 0.690521029504 | 1 |
| Cc1ccc(NC(=O)c2ccc(C#N)cn2)cc1S(N)(=O)=O | 0.688863375431 | 1 |
| Cc1ccc(NC(=O)c2ccncc2)cc1S(C)(=O)=O | 0.681247759053 | 2 |
| Cc1cc(NC(=O)c2ccccn2)ccc1S(N)(=O)=O | 0.680272108844 | 1 |
| Cc1ccc(NC(=O)c2ccc[nH]c2=O)cc1S(N)(=O)=O | 0.679851668727 | 5 |
| Cc1ccc(NC(=O)c2ccccn2)cc1S(=O)(=O)N(C)C | 0.679851668727 | 1 |
| Cc1ccc(NC(=O)c2cccnc2N)cc1S(=O)(=O)N(C)C | 0.669577874818 | 1 |
| Cc1ccc(NC(=O)c2ccc(N(C)C)nc2)cc1S(C)(=O)=O | 0.669506999391 | 1 |

αSMILES generated using both Paddy types.

**Table S2.** Paddy (Generational) Solutions using Tversky Similarity as Objective Function

| SMILES string | score | frequency |
|---|---|---|
| Cc1cnc(Nc2ccc(S(N)(=O)=O)cc2)nc1C | 0.775740479549 | 47 |
| Cc1cnc(Nc2ccc(S(N)(=O)=O)cc2)nc1N | 0.762300762301 | 1 |
| Cc1cnc(Nc2cccc(S(N)(=O)=O)c2)nc1C | 0.745542949757 | 6 |
| α Cc1ccc(NC(=O)c2cccnc2N)cc1S(N)(=O)=O | 0.738989062962 | 2 |
| Cc1cnc(Nc2cccc(S(N)(=O)=O)c2)nc1N | 0.7336523126 | 1 |
| α Cc1ccc(NC(=O)c2ccc(N)cn2)cc1S(N)(=O)=O | 0.730816077954 | 1 |
| Cc1cnc(Nc2ccccc2S(N)(=O)=O)nc1C | 0.727398683755 | 54 |
| Cc1ccc(NC(=O)c2ccccc2N)cc1S(N)(=O)=O | 0.724637681159 | 1 |
| Cc1cnc(Nc2ccc(S(C)(=O)=O)cc2)nc1C | 0.719969685487 | 1 |
| Cc1cnc(Nc2ccc(NS(C)(=O)=O)cc2)nc1C | 0.719969685487 | 1 |
| CCn1nccc1C(=O)Nc1ccc(C)c(S(N)(=O)=O)c1 | 0.717772035601 | 1 |
| Cc1cnc(Nc2ccccc2S(N)(=O)=O)nc1N | 0.715015321757 | 1 |
| Cc1cnc(Nc2ccc(S(=O)(=O)N(C)C)cc2)nc1 | 0.708661417323 | 1 |
| Cc1ncc(C)c(Nc2ccc(S(N)(=O)=O)cc2)n1 | 0.706582372629 | 1 |
| Cc1cc(NC(=O)c2nccn2C)ccc1S(N)(=O)=O | 0.700152207002 | 1 |
| Cc1cc(Nc2ccc(S(N)(=O)=O)cc2)nc(C)n1 | 0.694980694981 | 46 |
| Cc1ccccc1S(=O)(=O)Nc1ccccc1S(N)(=O)=O | 0.698080279232 | 3 |
| Cc1cnc(Nc2ccc(NS(C)(=O)=O)cc2)nc1 | 0.698080279232 | 2 |
| Cc1ccc(NC(=O)c2nccn2C)cc1S(C)(=O)=O | 0.690521029504 | 1 |
| Cc1ncnc(C)c1Nc1ccccc1S(N)(=O)=O | 0.682456844641 | 1 |
| Cc1nc(N)cc(Nc2ccc(S(N)(=O)=O)cc2)n1 | 0.681818181818 | 4 |
| Cc1ncc(Nc2ccc(S(N)(=O)=O)cc2)c(C)n1 | 0.681818181818 | 1 |
| Cc1ccccc1S(=O)(=O)Nc1ncc(N(C)C)cn1 | 0.681818181818 | 1 |
| Cc1ncc(C)c(Nc2ccccc2S(N)(=O)=O)n1 | 0.681247759053 | 1 |
| Cc1ccc(NC(=O)c2ccc(=O)n(C)n2)cc1S(C)(=O)=O | 0.679851668727 | 5 |

αSMILES generated using both Paddy types.

**Table S3.** Hyperopt Solutions using Tversky Similarity as Objective Function

| SMILES string | score | frequency |
|---|---|---|
| Cc1ccc(NC(=O)CC2CCCC2)cc1S(N)(=O)=O | 0.701530612245 | 249 |
| Cc1ccc(NC(=O)Cc2ccc(F)cc2)cc1S(N)(=O)=O | 0.701530612245 | 3 |
| Cc1ccc(NC(=O)Cc2ccc(Cl)cc2)cc1S(N)(=O)=O | 0.701530612245 | 2 |
| Cc1ccc(NC(=O)Cc2cccs2)cc1S(N)(=O)=O | 0.679851668727 | 29 |
| Cc1ccc(NC(=O)CCC2CCCC2)cc1S(N)(=O)=O | 0.679851668727 | 20 |
| Cc1ccc(NC(=O)Cc2ccsc2)cc1S(N)(=O)=O | 0.679851668727 | 10 |
| Cc1ccc(NC(=O)Cc2ccccc2F)cc1S(N)(=O)=O | 0.679851668727 | 7 |
| Cc1ccc(NC(=O)Cc2ccc[nH]2)cc1S(N)(=O)=O | 0.679851668727 | 2 |
| Cc1ccc(NC(=O)Cc2ccccc2Cl)cc1S(N)(=O)=O | 0.679851668727 | 2 |
| Cc1ccc(NC(=O)C2CC=CCC2)cc1S(N)(=O)=O | 0.679851668727 | 1 |
| Cc1ccc(NC(=O)CCSc2ccccn2)cc1S(N)(=O)=O | 0.678794461037 | 1 |
| Cc1ccc(NC(=O)Cc2cccc(F)c2)cc1S(N)(=O)=O | 0.669506999391 | 18 |
| Cc1ccc(NC(=O)Cc2cccc(Cl)c2)cc1S(N)(=O)=O | 0.669506999391 | 3 |
| Cc1ccc(NCc2cc(C#N)cs2)cc1S(N)(=O)=O | 0.669506999391 | 2 |

**Table S4.** Paddy (Population) Solutions using Custom Multi-Feature Objective Function

| SMILES string | score | frequency |
|---|---|---|
| Cc1ccc(NC(=O)c2cccc3ncccc23)cc1S(N)(=O)=O | 2.723916711 | 1 |
| Cc1ccc(NC(=O)c2cccc3ncccc23)cc1S(C)(=O)=O | 2.34155450409 | 1 |
| Cc1ccc(NC(=O)c2cc(C3CC3)[nH]n2)cc1S(N)(=O)=O | 2.12329223522 | 1 |
| α Cc1ccc(NC(=O)c2cccnc2N2CCCC2)cc1Cl | 2.11849336272 | 2 |
| α Cc1ccc(NC(=O)c2cccnc2N2CCCC2)cc1F | 2.10555596297 | 5 |
| α Cc1ccc(Nc2ncccc2C(=O)N2CCOCC2)cc1F | 2.0994818624 | 1 |
| Cc1ccc(NC(=O)c2ccccn2)cc1N1CCCC1=O | 2.01674640789 | 2 |
| α Cc1ccc(NC(=O)C2CC2)cc1Nc1ncccc1C#N | 2.01188802674 | 1 |
| Cc1ccc(NC(=O)c2cccc3ncccc23)cc1-n1cnnn1 | 2.00986738561 | 4 |
| Cc1ccc(S(C)(=O)=O)cc1NC(=O)c1cccc2ncccc12 | 2.0047763139 | 1 |
| Cc1ccc(NC(=O)c2ccnc(-n3ccnc3)c2)cc1Cl | 2.00087779538 | 1 |
| α Cc1ccc(NC(=O)c2cccnc2N2CCOCC2)cc1F | 1.99455379655 | 2 |
| Cc1ccc(NC(=O)c2ccccn2)cc1-n1cnnn1 | 1.99378763988 | 3 |
| Cc1ccc(Cl)cc1NC(=O)c1cccnc1N1CCCC1 | 1.99152386104 | 1 |
| α Cc1ccc(NC(=O)c2ccnc(-n3ccnc3)c2)cc1F | 1.99001940011 | 5 |
| α Cc1cc(NC(=O)c2cccnc2N2CCCC2)ccc1F | 1.97918828426 | 1 |
| Cc1cccc(NC(=O)Cn2cnc3c(cnn3C)c2=O)c1 | 1.97823445804 | 7 |
| α Cc1ccc(F)cc1NC(=O)c1cccnc1N1CCCC1 | 1.9762861328 | 1 |

αSMILES generated using both Paddy types.

**Table S5.** Paddy (Generational) Solutions using Custom Multi-Feature Objective Function

| SMILES string | score | frequency |
|---|---|---|
| Cc1ccc(Nc2ncccc2C(=O)NC2CC2)cc1F | 2.2654205205 | 40 |
| Cc1ccc(NC(=O)CSc2ncccc2N)cc1F | 2.24971825997 | 1 |
| Cc1ccc(Nc2ncccc2C(=O)N2CCCC2)cc1F | 2.22047936415 | 103 |
| Cc1ccc(NC(=O)C2CC2)cc1NCc1ccccn1 | 2.17662029154 | 1 |
| β Cc1ccc(NC(=O)NCc2cccnc2)cc1S(C)(=O)=O | 2.14919940829 | 1 |
| Cc1ccc(Nc2nc(C(=O)N3CCCC3)cs2)cc1F | 2.14446609203 | 2 |
| α Cc1ccc(NC(=O)c2cccnc2N2CCCC2)cc1Cl | 2.11849336272 | 4 |
| Cc1ccc(Nc2ncccc2C(=O)N2CCCCC2)cc1F | 2.1159438347 | 3 |
| Cc1cc(NC(=O)COc2cccc(F)c2)cc2ncccc12 | 2.11492395855 | 1 |
| Cc1cc(NC(=O)Cn2cccccc2=O)cc2ncccc12 | 2.11229538267 | 1 |
| Cc1ccc(Nc2ncccc2C(=O)N2CCOCC2)cc1Cl | 2.11082181245 | 1 |
| α Cc1ccc(NC(=O)c2cccnc2N2CCCC2)cc1F | 2.10555596297 | 313 |
| Cc1cc(Nc2ncccc2C(=O)N2CCCC2)ccc1F | 2.10255102342 | 3 |
| α Cc1ccc(Nc2ncccc2C(=O)N2CCOCC2)cc1F | 2.0994818624 | 11 |
| Cc1cccc(Nc2cc(C(=O)N3CCOCC3)ccn2)c1 | 2.07514478583 | 1 |
| Cc1ccc(NC(=O)NCCc2cccnc2)cc1S(C)(=O)=O | 2.05021322347 | 1 |
| Cc1ccc(NC(=O)c2ccnc(N3CCOCC3)c2)cc1F | 2.0491883323 | 5 |
| Cc1ccc(Nc2ncccc2C(=O)N2CCOCC2)cc1C | 2.04323215196 | 4 |
| Cc1ccc(NCc2cccnc2)cc1N1CCCC1=O | 2.04115602872 | 1 |
| Cc1ccc(Nc2ncccc2C(=O)NC2CCCC2)cc1F | 2.03727815308 | 2 |
| Cc1ccc(NC(=O)c2ccnc(-n3cncn3)c2)cc1Cl | 2.01691530646 | 1 |
| α Cc1ccc(NC(=O)C2CC2)cc1Nc1ncccc1C#N | 2.01188802674 | 2 |
| Cc1cccnc1CNC(=O)Nc1ccc2ncsc2c1 | 2.01038606454 | 1 |
| Cc1ccc(NC(=O)c2ccnc(-n3cncn3)c2)cc1F | 2.0060616961 | 5 |
| Cc1ccc(NC(=O)c2cccnc2N2CCOCC2)cc1Cl | 2.0057323405 | 1 |
| Cc1ccc(NC(=O)c2cc(N3CCOCC3)ccn2)cc1F | 1.99676506225 | 1 |
| Cc1cccc(Nc2ncccc2C(=O)N2CCOCC2)c1 | 1.99517654393 | 3 |
| α Cc1ccc(NC(=O)c2cccnc2N2CCOCC2)cc1F | 1.99455379655 | 43 |
| α Cc1ccc(NC(=O)c2ccnc(-n3ccnc3)c2)cc1F | 1.99001940011 | 13 |
| Cc1ccc(N)cc1NC(=O)c1cccnc1N1CCCC1 | 1.98938417372 | 1 |
| α Cc1cc(NC(=O)c2cccnc2N2CCCC2)ccc1F | 1.97918828426 | 2 |
| α Cc1ccc(F)cc1NC(=O)c1cccnc1N1CCCC1 | 1.9762861328 | 11 |
| Cc1ccc(NC(=O)c2cccnc2)cc1N1CCCC1=O | 1.97188960395 | 1 |

<sup>α</sup>SMILES generated using both Paddy types, <sup>β</sup>SMILES generated using both Paddy (generational) and Hyperopt.

**Table S6.** Hyperopt Solutions using Custom Multi-Feature Objective Function

| SMILES string | score | frequency |
|---|---|---|
| Cc1ccc(NC(=O)CSc2nccc(N)n2)cc1Cl | 2.35489428064 | 586 |
| Cc1ccc(Cl)cc1NC(=O)CSc1nccc(N)n1 | 2.25171293618 | 5 |
| Cc1ccc(NC(=O)CCn2cccn2)cc1S(C)(=O)=O | 2.19219435007 | 2 |
| Cc1ccc(OCC(=O)Nc2cccc3ncccc23)cc1C | 2.17151074779 | 3 |
| β Cc1ccc(NC(=O)NCc2cccnc2)cc1S(C)(=O)=O | 2.14919940829 | 1 |
| Cc1ccc(NC(=O)CC2CCCO2)cc1S(N)(=O)=O | 2.09813037529 | 1 |
| CC(=O)Nc1cccc(CNC(=O)c2ccc3c(C)ccnc3c2)c1 | 2.00645535494 | 1 |

βSMILES generated using both Paddy (generational) and Hyperopt.